\PassOptionsToPackage{table}{xcolor}
\documentclass{elsarticle}

\makeatletter
\myfooter[R]{\@date}
\makeatother

\usepackage[margin=1in]{geometry}

\usepackage{graphicx}
\usepackage{pstool}
\usepackage{wrapfig}
\usepackage{caption}
\usepackage{subcaption}

\usepackage[section]{placeins} 

\usepackage{fancyhdr} 
\usepackage{lastpage} 
\usepackage{extramarks} 

\usepackage{enumerate} 
\usepackage{paralist} 
\usepackage{amsmath, amsthm, amssymb, mathtools}
\usepackage{mathabx, pifont, stmaryrd} 
\usepackage[explicit]{titlesec} 
\usepackage{etoolbox} 
\usepackage{relsize}
\usepackage{bibentry} 
\makeatletter\let\saved@bibitem\@bibitem\makeatother 
\usepackage[colorlinks, bookmarksopen, bookmarksnumbered,
            citecolor=red,urlcolor=red]{hyperref} 
\makeatletter\let\@bibitem\saved@bibitem\makeatother 

\usepackage{algorithm}
\usepackage[noEnd=true,indLines=true]{algpseudocodex}

\usepackage{bm} 
\usepackage{amsfonts} 

\usepackage{tikz}
\usepackage{pgfplots}
\usepackage{pgfplotstable, filecontents, booktabs}
\pgfplotsset{compat=1.9}

\usetikzlibrary{pgfplots.groupplots}
\usepgfplotslibrary{fillbetween}
\usetikzlibrary{calc,fit,matrix,arrows,automata,positioning,shapes}
\usetikzlibrary{arrows.meta}

\pgfplotsset{select coords between index/.style 2 args={
    x filter/.code={
        \ifnum\coordindex<#1\fi
        \ifnum\coordindex>#2\fi
    }
}}

\tikzset{
 invisible/.style={opacity=0},
 visible on/.style={alt={#1{}{invisible}}},
 alt/.code args={<#1>#2#3}{%
   \alt<#1>{\pgfkeysalso{#2}}{\pgfkeysalso{#3}}
 },
}

\usepackage{setspace}
\usepackage{csquotes}
\usepackage[normalem]{ulem}

\newbool{fastcompile}
\setbool{fastcompile}{false}

\newcommand{\EE}{{\mathbb E}}
\newcommand{\RR}{{\mathbb R}}

\newcommand{\cM}{{\mathcal M}}

\newcommand{\cU}{{\mathcal U}}

\newcommand{\ctrl}{\alpha}
\newcommand{\ctrldim}{k}
\newcommand{\dom}{\mathcal{D}}

\usepackage{empheq}
\usepackage{caption} 
\usepackage{dsfont} 

  \definecolor{Green}{RGB}{198, 239, 206}
  \definecolor{Salmon}{RGB}{255, 199, 206}
  \definecolor{SkyBlue}{RGB}{189, 215, 238}

\begin{document}

\title{\texorpdfstring{Picard-Based Acceleration of Newton Continuation for Mean Field Game PDE Systems}{Picard-Based Acceleration of Newton Continuation for Mean Field Game PDE Systems}}

\author[rvt1,rvt2]{Mathieu Laurière\fnref{fn1}}
\ead{mathieu.lauriere@nyu.edu}

\author[rvt2]{Andrew Shi\fnref{fn2}\corref{cor1}}
\ead{andrewshi@math.berkeley.edu}

\address[rvt1]{NYU Shanghai Frontiers Science Center of Artificial Intelligence and Deep Learning}
\address[rvt2]{NYU-ECNU Institute of Mathematical
Sciences at NYU Shanghai}
\cortext[cor1]{Corresponding author}

\fntext[fn1]{Assistant Professor, Department of Mathematics and Department of Data Science, New York University Shanghai}
\fntext[fn2]{Postdoctoral Instructor, Department of Mathematics, New York University Shanghai}


\begin{abstract}
We develop a method that uses Picard iterations to accelerate Newton continuation for the semi-implicit finite-difference discretization of forward-backward partial differential equation (PDE) systems arising in mean field games (MFGs).
We first investigate the computational properties of the Picard and Newton methods, which are widely used separately in the MFG literature but whose comparative cost and robustness across different regimes remain insufficiently explored and documented.
The Picard method uses an outer fixed-point iteration that alternates a forward Fokker-Planck solve and a backward Hamilton-Jacobi-Bellman solve.
The Newton method instead applies Newton's method directly to the coupled nonlinear space-time system.
Across one- and two-dimensional MFG benchmarks, Picard offers substantial computational savings in favorable regimes, but may fail at sufficiently low viscosity or require strong damping under temporal shocks, making Newton continuation preferable.
With parameter continuation in the viscosity parameter, the Newton method is more robust in these regimes, at the cost of larger coupled linear systems.
We relate these trade-offs to the residuals, Jacobian blocks, and sparsity structures produced by separable, local nonseparable, and nonlocal Hamiltonians.
We then demonstrate how to combine inexpensive Picard iterations with Newton continuation in a hybrid method to reduce the total computational cost for a two-dimensional double-well problem.
\end{abstract}

\maketitle

\section{Introduction}
\label{sec:intro}

Mean field game (MFG) theory provides a framework for differential games involving a continuum of indistinguishable agents whose individual strategies are coupled through the population distribution.
We consider three classes of MFG systems: finite-horizon, infinite-horizon discounted, and ergodic. Each is characterized by a coupled system of partial differential equations (PDEs).
In the finite-horizon case, the system has a forward--backward structure: the backward-in-time Hamilton--Jacobi--Bellman (HJB) equation determines the optimal value function for a representative agent subject to a terminal cost, while the forward Fokker--Planck (FP), or Kolmogorov, equation describes the evolution of the population distribution from an initial configuration.

Mean field games were introduced independently by Lasry and Lions~\cite{Lasry2006MFG1,Lasry2006MFG2,Lasry2007MFGJPN} and by Huang, Malhamé, and Caines~\cite{Huang2006MFG}.
For general introductions to MFG theory, see~\cite{Gueant2011ParisPrinceton,bensoussan2013mfgmfc,Cardaliaguet2020MFGTheory}; for a systematic treatment of the probabilistic approach, see~\cite{CarmonaDelarue2018ProbabilisticI}.
Applications include economic competition~\cite{Chan2014BertrandCournot}, crowd dynamics~\cite{achdou2018crowdmotion}, and segregation models~\cite{achdou2017segregation}.

The finite-difference method of Achdou and Capuzzo-Dolcetta~\cite{achdou2010numerical} is the starting point for the present work.
Convergence results for this discretization were established in~\cite{achdou2013FDconvergence,achdou2016FDconvergence}; related developments include its extension to the planning problem~\cite{achdou2012planning}.
For variational formulations, augmented-Lagrangian methods were studied in~\cite{Benamou2015augmentedlagrangian,Andreev2017preconditioning}, primal-dual algorithms based on the Chambolle--Pock method in~\cite{Chambolle2010ChambollePock,BricenoArias2019Chambelle}, and proximal methods in~\cite{BricenoArias2018proximal}.
Other discretizations include semi-Lagrangian and related Lagrange--Galerkin schemes~\cite{carlini2014semilagrangianfirstorder,Carlini2024LagrangeGalerkin,Carlini2026monotone}, finite-element methods~\cite{Fu2023highorderFEM,Osborne2025FEMNondifferentiableHamiltonians}, and mesh-free methods based on Gaussian processes and Fourier features~\cite{Mou2022fourier}.
Machine-learning methods have also been developed for high-dimensional MFG and mean field control problems~\cite{Ruthotto2020MLhighdim,Carmona2021Ergodic,Lin2021StochasticMFG}.
For broader reviews of numerical methods for MFGs, see~\cite{achdou2020numericalaspects,lauriere2021AMSnumerical,hu2024recent}.

Full space-time Newton solvers for discrete MFG systems appear in~\cite{achdou2012planning,achdou2012preconditioning}, while alternating forward--backward fixed-point iterations related to the Picard method considered here appear in~\cite{achdou2020numericalaspects,Li2021multiscale}.
Li, Fan, and Ying~\cite{Li2021multiscale} establish a local spectral-radius condition for their alternating-sweeping method and compare its runtime with Newton on a one-dimensional benchmark.
Carlini and Zorkot~\cite{Carlini2027Newton} compare Newton linearization before and after discretization using finite-difference and semi-Lagrangian schemes, and use an Armijo line search to globalize the iterations.
The solver structures themselves are not new. Our main contributions are:
\begin{enumerate}
    \item a comparison of Picard and Newton solvers for the same semi-implicit finite-difference discretization across one- and two-dimensional MFG benchmarks, identifying the complementary regimes that motivate the hybrid strategy: Picard can provide inexpensive initial solves, while Newton continuation supplies the robustness needed to reach lower viscosities;
    \item an implementation-level account of the residuals, Jacobian blocks, and sparsity structures that determine solver cost for separable, local nonseparable, and nonlocal Hamiltonians;
    \item a hybrid continuation method that uses inexpensive Picard iterations to reduce the number of coupled Newton solves and the total computational cost.
\end{enumerate}

The rest of the paper is organized as follows.
Section~\ref{sec:eqnsanddisc} recalls the finite-horizon, discounted, and ergodic MFG systems used in the benchmarks.
Section~\ref{sec:numericalmethods} presents the finite-difference discretization, the Newton and Picard solvers, the continuation strategy, and the hybrid protocol.
Section~\ref{sec:numexp} first compares the methods on validation tests and identifies regimes in which one method is clearly preferable, making a hybrid approach unnecessary or impractical. It then evaluates the computational savings obtained by combining the methods in the hybrid strategy.
Section~\ref{sec:conc} summarizes the computational gains and limitations of the hybrid method and discusses directions toward scalable MFG solvers, including matrix-free implementations and intermediate Picard/Newton variants.

\section{Governing equations}
\label{sec:eqnsanddisc}

\subsection{MFG PDE systems for different optimal control formulations}
\label{sec:eqnsanddisc:MFGPDE}

We recall the three MFG PDE systems used in the numerical experiments: the finite-horizon, infinite-horizon discounted, and ergodic systems.
Let $\dom$ denote the state domain, let $d$ be its dimension, and let $k$ be the control dimension.
For the local models considered below, the population argument in $H(x,m,p)$ is the pointwise density; the nonlocal dependence used later is introduced separately in Section~\ref{sec:numexp:nonlocalcongestion}. Given a running cost $f$, terminal cost $g$, and drift $b$, the Hamiltonian is
\begin{equation}  \label{eq:Hamiltonian-from-Lagrangian}
    H(x,m,p) = \sup_{\ctrl \in \RR^\ctrldim} \left\{ - \langle b(x,m,\ctrl), p \rangle - f(x, m, \ctrl) \right\}.
\end{equation}
We restrict attention to data for which the supremum in \eqref{eq:Hamiltonian-from-Lagrangian} is attained by a unique feedback control and the derivatives of $H$ required by the Newton linearizations below are well defined. We denote the unique maximizer by $\hat{\ctrl}(x,m,p)$. Each system below is understood to be supplemented with the spatial boundary conditions specified for the corresponding numerical experiment in Section~\ref{sec:numexp}.

\subsubsection{Finite-horizon optimal control}
\label{sec:eqnsanddisc:MFGPDE:finitehorizon}

The finite-horizon system was formally derived in \cite{Lasry2006MFG2} and analyzed extensively in \cite{Cardaliaguet2010NotesMFG}. 
For simplicity, let the initial time be $0$, and let $T > 0$ be the fixed time horizon. Given a fixed flow of measures $m = (m_t)_{t \in [0, T]}$ admitting densities, the representative agent minimizes
\begin{equation}
    J^{MFG}_{m}(\ctrl) = \EE \left[ \int_0^T f(X_t, m_t(X_t), \ctrl_t) \, dt + g(X_T, m_T) \right],
\end{equation}
subject to the dynamics
\begin{equation}
    d X_t = b(X_t, m_t(X_t), \ctrl_t) \, dt + \sigma \, d W_t, \quad X_0 \sim m_0.
\end{equation}
Here, $W$ is a $d$-dimensional standard Brownian motion and $\sigma$ is the volatility coefficient, assumed to be scalar and positive for simplicity. 
When a population measure admits a density, we use the same notation $m_t$ for the measure and its density.
The same identification extends to the discrete setting in Section~\ref{sec:numericalmethods:FD-discretization}, where $M_{i,j}$ denotes the density value at the grid point $x_{i,j}$.
For the nonlocal congestion model in
Section~\ref{sec:numexp:nonlocalcongestion}, the pointwise density
$m_t(X_t)$ is replaced by the averaged density $(K\ast m_t)(X_t)$.

The resulting system of coupled PDEs is as follows:
\begin{equation} \label{eq:forward-backward-PDE}
\begin{cases}
    -\frac{\partial u}{\partial t}(t, x) - \nu \Delta u(t, x) + H(x, m(t, x), \nabla u(t, x)) = 0, & (t,x) \in (0,T) \times \dom, \\
    \frac{\partial m}{\partial t}(t, x) - \nu \Delta m(t, x) - \text{div}\left(m(t, x) D_p H(x, m(t, x), \nabla u(t, x)) \right) = 0, & (t,x) \in (0,T) \times \dom, \\
    m(0,x) = m_0(x), \quad u(T,x) = g(x, m(T)), & x \in \dom
\end{cases}
\end{equation}
where $\nu = \frac{\sigma^2}{2}$.
By the envelope theorem, the term $-D_p H(x, m, \nabla u)$ equals the optimal drift $b(x, m, \hat{\ctrl})$. 
The FP equation describes the evolution of the density under this optimal drift.

\subsubsection{Infinite-horizon discounted MFG system}
\label{sec:eqnsanddisc:MFGPDE:discounted}

We next consider the stationary MFG system associated with an infinite-horizon discounted control problem, with discount factor $\rho>0$.

For a fixed stationary population density $m$ and initial state $X_0=x$, the representative agent minimizes the infinite-horizon discounted cost
\begin{equation}
    J^{MFG}_{m}(\ctrl) = \EE \left[ \int_0^\infty e^{-\rho t} f(X_t,m(X_t),\ctrl_t)\,dt \right],
\end{equation}
where $X$ satisfies $dX_t=b(X_t,m(X_t),\ctrl_t)\,dt+\sigma\,dW_t$ and $m$ is held fixed. The stationary equilibrium is characterized by the following coupled PDE system:
\begin{equation} \label{eq:stationary-discounted-PDE-system}
\begin{cases}
    \rho u(x) - \nu\Delta u(x) + H(x, m(x), \nabla u(x)) = 0, & x \in \dom, \\
    -\nu\Delta m(x) - \text{div}\left( m(x)D_p H(x, m(x), \nabla u(x)) \right) = 0, & x \in \dom, \\
    \int_{\dom} m(x) \, dx = 1.
\end{cases}
\end{equation}

The stationary system \eqref{eq:stationary-discounted-PDE-system} characterizes the discounted equilibrium and admits closed-form solutions in special cases (Section~\ref{sec:numexp:discounted}), making it a useful analytical benchmark. 
For numerical computation, however, we work with a finite-horizon forward-backward formulation that reuses the time-marching machinery developed for the finite-horizon system \eqref{eq:forward-backward-PDE}. 
Introducing a computational horizon $T$ and an arbitrary terminal cost $g$, we obtain
\begin{equation} \label{eq:forward-backward-discounted-PDE-system}
\begin{cases}
    -\frac{\partial u}{\partial t}(t, x) +\rho u(t, x) -\nu \Delta u(t, x) + H(x, m(t, x), \nabla u(t, x)) = 0, & (t,x) \in (0,T) \times \dom, \\
    \frac{\partial m}{\partial t}(t, x) - \nu \Delta m(t, x) - \text{div}\left(m(t, x) D_p H(x, m(t, x), \nabla u(t, x)) \right) = 0, & (t,x) \in (0,T) \times \dom, \\
    m(0,x) = m_0(x), \quad u(T,x) = g(x, m(T)), & x \in \dom.
\end{cases}
\end{equation}

\subsubsection{Infinite-horizon ergodic MFG system}
\label{sec:eqnsanddisc:MFGPDE:ergodic}

The ergodic formulation was introduced in \cite{Lasry2006MFG1} and is further detailed in \cite{Gueant2011ParisPrinceton}. 
Under assumptions ensuring ergodicity of the optimally controlled diffusion, the long-run average cost is independent of the initial state.

For a fixed stationary population density $m$, the representative agent minimizes the long-run average cost
\begin{equation}
    J^{MFG}_{m}(\ctrl) = \limsup_{T \to \infty} \frac{1}{T} \EE \left[ \int_0^T f(X_t, m(X_t), \ctrl_t) \, dt \right].
\end{equation}
Here, $X$ satisfies $dX_t=b(X_t,m(X_t),\ctrl_t)\,dt+\sigma\,dW_t$. Let $\lambda=\inf_{\ctrl}J^{MFG}_{m}(\ctrl)$ denote the corresponding optimal long-run average cost. At an MFG equilibrium, $m$ is invariant under an optimal feedback control.
Unless it is fixed by boundary data, the HJB equation determines $u$ only up to an additive constant. In \eqref{eq:stationary-ergodic-PDE-system}, we select the representative satisfying $\int_{\dom}u(x)\,dx=0$; in the Dirichlet validation problem of Section~\ref{sec:numexp:ergodic}, the boundary data fix this constant instead. 
The resulting equilibrium triplet $(\lambda,u,m)$ satisfies the following coupled PDE system:
\begin{equation} \label{eq:stationary-ergodic-PDE-system}
\begin{cases}
    \lambda - \nu \Delta u(x) + H(x, m(x), \nabla u(x)) = 0, & x \in \dom, \\
    -\nu \Delta m(x) - \text{div}\left( m(x) D_p H(x, m(x), \nabla u(x)) \right) = 0, & x \in \dom, \\
    \int_{\dom} m(x) \, dx = 1, \quad \int_{\dom} u(x) \, dx = 0.
\end{cases}
\end{equation}

These formulations play different roles below. The finite-horizon systems \eqref{eq:forward-backward-PDE} and \eqref{eq:forward-backward-discounted-PDE-system} are the main time-dependent problems solved by the Newton and Picard methods. The stationary discounted and ergodic systems provide analytical or discrete reference solutions for the validation experiments in Sections~\ref{sec:numexp:discounted} and~\ref{sec:numexp:ergodic}.

\subsection{Connections between the formulations}
\label{sec:eqnsanddisc:connections}

The discounted, ergodic, and finite-horizon systems are used in the numerical experiments through two standard asymptotic links.
To compare the HJB components, fix a density $m$ and write
\[
\mathcal{L}_m(\phi)=-\nu\Delta \phi+H(x,m(x),\nabla\phi).
\]
The stationary discounted and ergodic HJB equations then take the forms
\[
\rho u+\mathcal{L}_m(u)=0,
\qquad
\lambda+\mathcal{L}_m(u)=0.
\]
For the coupled MFG systems, under assumptions ensuring the vanishing-discount limit, the discounted value function $u_\rho$ satisfies $\rho u_\rho\to\lambda$, and the normalized family $u_\rho-\lambda/\rho$ remains bounded; see, e.g., \cite{CirantPorretta2021Turnpike}.
In Section~\ref{sec:numexp:ergodic}, we test a discrete version of this relation on a fixed truncated domain and grid.

Under suitable structural assumptions, the finite-horizon and ergodic systems exhibit a turnpike property: for large horizons, the finite-horizon trajectory remains close to the stationary ergodic state away from the initial and terminal boundary layers \cite{CirantPorretta2021Turnpike}. 
At the level of the value function, this corresponds formally to
\begin{equation} \label{eq:ansatz-for-finite-horizon}
u(t,x)\approx \lambda(T-t)+\bar u(x),
\end{equation}
where $(\lambda,\bar u,\bar m)$ solves the corresponding ergodic problem.
We use this relation only as a numerical validation device: when no finite-horizon closed form is available, midpoint slices of long-horizon computations can be compared with stationary references.

\section{Numerical methods}
\label{sec:numericalmethods}
 
\subsection{Finite difference discretization}
\label{sec:numericalmethods:FD-discretization}

We present the discretization in two spatial dimensions following \cite{achdou2010numerical}; the one-dimensional case follows by restriction. 
Let $N_h$ be a positive integer, let $h$ denote the mesh spacing, and let $\dom_h^d$ denote the resulting uniform grid on $\dom$ with $N_h$ points in each direction.
In the one-dimensional Dirichlet benchmarks, however, $N_h$ denotes the number of spatial intervals, so the vectors and matrices have $N_h+1$ nodal entries; for the periodic and cell-centered benchmarks, $N_h$ denotes the number of grid points in each coordinate direction. For subsequent dimension counts, let $N_{\mathrm{sp}}$ denote the number of
spatial degrees of freedom: $N_{\mathrm{sp}}=N_h+1$ in the one-dimensional
Dirichlet cases and $N_{\mathrm{sp}}=N_h^d$ otherwise.
For simplicity, we take the domains to be squares in 2D and discretize them with the same mesh size in both directions.
Let $x_{i, j}$ denote a point in $\dom_h^2$; the values of $u$ and $m$ at $x_{i, j}$ are denoted by $U_{i, j}$ and $M_{i, j}$, respectively. 
For temporal discretization, let $N_T$ be a positive integer, set $\Delta t = T/N_T$, and define $t_n = n \Delta t$ for $n = 0, \dots, N_T$. 
The values of $u$ and $m$ at $(x_{i, j}, t_n)$ are denoted by $U_{i, j}^n$ and $M_{i, j}^n$, respectively. 

To approximate the spatial gradient $\nabla u$, we rely on a four-point stencil combining forward and backward differences in both spatial dimensions. 
We evaluate the numerical Hamiltonian $\tilde{H}$ by assigning the components of the discrete gradient $[D_hW]_{i,j}$ to the arguments $q_1, \dots, q_4$. 
Specifically, $q_1$ and $q_2$ correspond to the forward and backward differences in the $x_1$-direction, while $q_3$ and $q_4$ correspond to the forward and backward differences in the $x_2$-direction.

We define the following operators:
\begin{subequations}
\label{eq:finite_differences}
\begin{align}
(D_t W)^n &= \frac{1}{\Delta t}(W^{n+1} - W^n) \\
(D_1^+ W)_{i,j} & = \frac{W_{i+1,j} - W_{i,j}}{h} \quad \text{and} \quad (D_2^+ W)_{i,j} = \frac{W_{i,j+1} - W_{i,j}}{h} \\
(\Delta_h W)_{i,j} &= -\frac{1}{h^2}(4W_{i,j} - W_{i+1,j} - W_{i-1,j} - W_{i,j+1} - W_{i,j-1}) \\
[D_h W]_{i,j} &= \big( (D_1^+ W)_{i,j}, (D_1^+ W)_{i-1,j}, (D_2^+ W)_{i,j}, (D_2^+ W)_{i,j-1} \big)^T.
\end{align}
\end{subequations}

The discrete Hamiltonian $\tilde{H}: \dom \times \RR \times \RR^4 \rightarrow \RR$ must satisfy the following standard assumptions for monotone schemes \cite{achdou2010numerical, achdou2013FDconvergence}:

\begin{itemize}
\item[$(\textbf{H}_1)$] \textbf{Monotonicity:} $\tilde{H}$ is nonincreasing with respect to $q_1$ and $q_3$ and nondecreasing with respect to $q_2$ and $q_4$.
\item[$(\textbf{H}_2)$] \textbf{Consistency:} $\tilde{H}(x, m, p_1, p_1, p_2, p_2) = H(x, m, (p_1, p_2))$ for all $x \in \dom$, $m \in \mathbb{R}$, and $(p_1, p_2) \in \mathbb{R}^2$.
\item[$(\textbf{H}_3)$] \textbf{Regularity:} $\tilde{H}$ is at least $\mathcal{C}^2$ with respect to $m$ and $q$.
\item[$(\textbf{H}_4)$] \textbf{Convexity:} The mapping $(q_1, q_2, q_3, q_4) \mapsto \tilde{H}(x, m, q_1, q_2, q_3, q_4)$ is convex.
\end{itemize}

In this work, we use the Godunov scheme. 
For the Hamiltonians considered below, which depend on the spatial gradient only through its magnitude and are nondecreasing in that magnitude (i.e., $H(x, m, p) = f(x, m, |p|)$ with $f$ nondecreasing in $|p|$), the Godunov scheme simplifies to the standard isotropic upwind formulation:
\begin{equation} \label{eq:godunov-hamiltonian}
    \tilde{H}(x, m, q_1, q_2, q_3, q_4) = H\left(x, m,\sqrt{(q_1^-)^2 + (q_3^-)^2 + (q_2^+)^2 + (q_4^+)^2}\right)    
\end{equation}
where $(\cdot)^+$ and $(\cdot)^-$ denote the projections onto $\mathbb{R}_+$ and $\mathbb{R}_-$, respectively. 
Other suitable numerical Hamiltonians, such as the Lax--Friedrichs Hamiltonian presented in Example 2 of \cite{achdou2010numerical}, could also be used. 

Therefore, the discrete HJB equation corresponding to \eqref{eq:forward-backward-PDE} is

\begin{equation} \label{eq:discrete-HJB-finitediff}
\begin{cases}
    - (D_t U_{i, j})^{n} - \nu (\Delta_h U^{n})_{i, j}
	+ \tilde{H}(x_{i, j}, M^{n+1}_{i, j}, [D_h U^n]_{i, j}) 
    = 0 , \\ 
	U^{N_T}_{i, j} = g(x_{i, j}, M^{N_T}_{i, j}).
\end{cases}
\end{equation}

For the discounted problem, the discrete HJB equation contains an additional term:

\begin{equation} \label{eq:discrete-HJB-finitediff-discounted}
\begin{cases}
    - (D_t U_{i, j})^{n} - \nu (\Delta_h U^{n})_{i, j}
	+ \rho U_{i, j}^n+\tilde{H}(x_{i, j}, M^{n+1}_{i, j}, [D_h U^n]_{i, j}) 
    = 0 , \\ 
	U^{N_T}_{i, j} = g(x_{i, j}, M^{N_T}_{i, j}).
\end{cases}
\end{equation}

Boundary conditions are imposed separately for each numerical experiment.

We discretize the FP equation in \eqref{eq:forward-backward-PDE} through its weak form. 
Integration by parts gives
\begin{equation}
-\int_{\dom}\text{div}\left(m\frac{\partial H}{\partial p}(x, m, \nabla u)\right)w = \int_{\dom}m\frac{\partial H}{\partial p}(x, m, \nabla u) \cdot \nabla w
\end{equation}
where the boundary term $\int_{\partial D} (m D_p H(x, m, \nabla u)) \cdot n \, w \, dS$ vanishes under any of the standard boundary conditions considered in this work: periodic (the boundary $\partial D$ is empty on the torus), homogeneous Neumann on the FP flux ($m D_p H \cdot n = 0$), or homogeneous Dirichlet on $m$. 
We approximate this expression by
\begin{equation}
h^2\sum_{i, j}M_{i, j}\nabla_q\tilde{H}(x_{i, j}, M_{i, j}, [D_hU]_{i, j}) \cdot [D_hW]_{i, j}.
\end{equation}
Adopting the notation
\begin{equation} \label{eq:transport-operator}
\begin{aligned}
\mathcal{B}_{i, j}(U, M) &= \frac{1}{h} \Bigg[ M_{i, j}\frac{\partial \tilde{H}}{\partial q_1}(x_{i, j}, M_{i, j},[D_hU]_{i, j}) - M_{i-1, j}\frac{\partial \tilde{H}}{\partial q_1}(x_{i-1, j}, M_{i-1, j}, [D_hU]_{i-1, j}) \\
&\quad + M_{i+1, j}\frac{\partial \tilde{H}}{\partial q_2}(x_{i+1, j}, M_{i+1, j}, [D_hU]_{i+1, j}) - M_{i, j}\frac{\partial \tilde{H}}{\partial q_2}(x_{i, j}, M_{i, j}, [D_hU]_{i, j}) \\
&\quad + M_{i, j}\frac{\partial \tilde{H}}{\partial q_3}(x_{i, j}, M_{i, j}, [D_hU]_{i, j}) - M_{i, j-1}\frac{\partial \tilde{H}}{\partial q_3}(x_{i, j-1}, M_{i, j-1}, [D_hU]_{i, j-1}) \\
&\quad + M_{i, j+1}\frac{\partial \tilde{H}}{\partial q_4}(x_{i, j+1}, M_{i, j+1}, [D_hU]_{i, j+1}) - M_{i, j}\frac{\partial \tilde{H}}{\partial q_4}(x_{i, j}, M_{i, j}, [D_hU]_{i, j}) \Bigg].
\end{aligned}
\end{equation}

The eight terms in $\mathcal{B}_{i,j}$ form the discrete divergence of $m D_p \tilde{H}$: each of the four coordinate-direction pairs contributes a flux difference across one face of the cell at $(i,j)$, with the $q_\ell$ index selecting the one-sided difference of $u$ consistent with that face.
This yields the following discretization of the FP equation in \eqref{eq:forward-backward-PDE}:

\begin{equation} \label{eq:discrete-FP-finitediff}
\begin{cases}
   (D_t M_{i, j})^{n} - \nu (\Delta_h M^{n+1})_{i, j} - \mathcal{B}_{i, j}(U^{n}, M^{n+1}) = 0 \\
	M^{0}_{i, j} = m_0(x_{i, j}).
\end{cases}
\end{equation}

Together, these equations yield the discrete forward-backward MFG system corresponding to \eqref{eq:forward-backward-PDE}.
Given $M^0$ and $U^{N_T}$, find $(U^n,M^n)$ for $0\leq n\leq N_T$ such that, for all $0\leq n<N_T$ and $0\leq i,j<N_h$,

\begin{subequations}\label{eq:MFG-finitedifference-fb}
     \begin{empheq}{align}
     \label{eq:MFG-finitedifference-fb:HJB}     
-\frac{U_{i, j}^{n+1} - U_{i, j}^n}{\Delta t} - \nu
(\Delta_hU^{n})_{i, j}+\tilde{H}(x_{i, j}, M_{i, j}^{n+1}, [D_hU^{n}]_{i, j}) &= 0 \\ 
\label{eq:MFG-finitedifference-fb:FP}     
\frac{M_{i, j}^{n+1} - M_{i, j}^n}{\Delta t} - \nu
(\Delta_hM^{n+1})_{i, j}-\mathcal{B}_{i, j}(U^{n}, M^{n+1}) &= 0. 
     \end{empheq}
\end{subequations}
When solving \eqref{eq:MFG-finitedifference-fb:HJB} for $U^n$, the mass $M^{n+1}$ is treated as known; when solving \eqref{eq:MFG-finitedifference-fb:FP} for $M^{n+1}$, the value $U^n$ is treated as known. 
Each equation is thus implicit in its own unknown but explicit in the coupling; we therefore refer to this treatment as a semi-implicit discretization.

For the discounted case, \eqref{eq:MFG-finitedifference-fb:HJB} becomes
\begin{equation}
-\frac{U_{i,j}^{n+1} - U_{i,j}^{n}}{\Delta t} - \nu (\Delta_h U^n)_{i,j} + \rho U_{i,j}^n + \tilde{H}(x_{i,j}, M_{i,j}^{n+1}, [D_h U^n]_{i,j}) = 0.
\end{equation}

\subsection{Newton and Picard methods for the forward-backward problem}
\label{sec:numericalmethods:Newton}

We now consider two strategies for solving the discrete forward-backward MFG system \eqref{eq:MFG-finitedifference-fb}. 

\begin{enumerate}
    \item The \textbf{Newton method} applies Newton's method simultaneously to all unknowns $U$ and $M$ in the full system \eqref{eq:MFG-finitedifference-fb}. This captures the full coupling between the HJB and FP equations and results in a significantly larger linear system.
    \item The \textbf{Picard method} decouples $U$ and $M$ using an outer Picard iteration. 
    For a fixed $U$, it marches a linearized FP equation forward from $m_0$, with one linear solve per timestep; for non-separable Hamiltonians, density-dependent transport coefficients are evaluated at the preceding Picard iterate, as described below.
    For a fixed $M$, it solves the nonlinear HJB equation \eqref{eq:MFG-finitedifference-fb:HJB} for the value function $U$ by starting with the prescribed terminal condition $g_T$ and marching backward in time; in our implementation, Newton's method is used at each timestep to handle the nonlinear Hamiltonian.
\end{enumerate}

For the Newton method, we stack the residuals of \eqref{eq:MFG-finitedifference-fb:HJB} and \eqref{eq:MFG-finitedifference-fb:FP} over all space-time grid points into a single vector-valued function
\begin{equation}
\varphi(U, M) = (\varphi_U(U, M), \varphi_M(U, M))^\top,
\end{equation}
where $\varphi_U(U, M) \in \mathbb{R}^{(N_T+1)N_{\mathrm{sp}}}$ collects the HJB residuals and $\varphi_M(U, M) \in \mathbb{R}^{(N_T+1)N_{\mathrm{sp}}}$ collects the FP residuals. A zero of $\varphi$ is equivalent to a solution of the discrete system \eqref{eq:MFG-finitedifference-fb} together with its initial and terminal conditions.

Let $J_\varphi$ denote the Jacobian of $\varphi$. Newton's method starts from an initial guess $(U^{(0)}, M^{(0)})$ and at each iterate $k$ solves
\begin{equation*}
J_\varphi(U^{(k)}, M^{(k)}) \, \delta W^{(k)} = -\varphi(U^{(k)}, M^{(k)})
\end{equation*}
for the update $\delta W^{(k)} = (\delta U^{(k)}, \delta M^{(k)})^\top$, then sets $(U^{(k+1)}, M^{(k+1)})^\top = (U^{(k)}, M^{(k)})^\top + \delta W^{(k)}$.

Each step amounts to solving a linear system of the form
\begin{equation}
\begin{pmatrix} A_{\cU,\cU} & A_{\cU,\cM} \\ A_{\cM,\cU} & A_{\cM,\cM} \end{pmatrix}
\begin{pmatrix} \delta U \\ \delta M \end{pmatrix}
=
-\begin{pmatrix} \varphi_U \\ \varphi_M \end{pmatrix},
\end{equation}
where the Jacobian blocks are evaluated at the current iterate $(U^{(k)},M^{(k)})$ and are given by
\begin{align*}
A_{\cU,\cM}(U, M) &= \nabla_M \varphi_U(U, M), & A_{\cU,\cU}(U, M) &= \nabla_U \varphi_U(U, M), \\
A_{\cM,\cU}(U, M) &= \nabla_U \varphi_M(U, M), & A_{\cM,\cM}(U, M) &= \nabla_M \varphi_M(U, M).
\end{align*}

We describe these four matrices in some detail. 

\begin{equation} \label{eq:globalNewton-Jacobian_AUU}
	A_{\cU, \cU}
	=
\begin{pmatrix} 
D_0 & -\frac{1}{\Delta t} \text{Id}_{N_{\mathrm{sp}}} & 0 & \dots & 0 \\
0 & D_1 & -\frac{1}{\Delta t} \text{Id}_{N_{\mathrm{sp}}} & \ddots & \vdots \\
\vdots & \ddots & \ddots & \ddots & 0 \\
0 & \dots & 0 & D_{N_T-1} & -\frac{1}{\Delta t} \text{Id}_{N_{\mathrm{sp}}} \\
0 & \dots & 0 & 0 & \text{Id}_{N_{\mathrm{sp}}}
\end{pmatrix}
\end{equation}

where, in two dimensions, $D_n$ represents the discrete operator

\begin{equation}
	Z = (Z_{i,j})_{i,j}
	\mapsto
	\left( \frac{1}{\Delta t} Z_{i,j} - \nu (\Delta_h Z)_{i,j} + [D_h Z]_{i,j} \cdot \nabla_q \tilde{H}(x_{i,j}, M_{i,j}^{n+1}, [D_h U^{(\mathtt{k}), n}]_{i,j})
	\right)_{i,j},
\end{equation}

which comes from the linearization of the discrete HJB equation~\eqref{eq:MFG-finitedifference-fb:HJB}. 
The identity matrix on the bottom right enforces the terminal condition $U^{N_T} = g(x, M^{N_T})$.\footnote{Our block layout differs from related references due to indexing and boundary-condition conventions. In equation (72) of \cite{achdou2012planning} the identity block sits at the top left and the $-\frac{1}{\Delta t}\text{Id}_{N_{\mathrm{sp}}}$ blocks on the first subdiagonal, because \cite{achdou2012planning} poses the forward-backward system with an initial condition on $u$ and a terminal condition on $m$. In equation (15) of \cite{achdou2012preconditioning} the identity block is absent from $A_{\cU, \cU}$ because BC enforcement is handled outside the Newton iteration.}

For the locally separable Hamiltonian $H(x,m,p)=H_0(x,p)-f_0(x,m)$, linearization preserves the adjoint relationship between the HJB and FP equations, yielding the following form of $A_{\cM,\cM}$:

\begin{equation}  \label{eq:globalNewton-Jacobian_AMM}
	A_{\cM, \cM}
	=
\begin{pmatrix} 
\text{Id}_{N_{\mathrm{sp}}} & 0 & \dots & 0 & 0 \\
-\frac{1}{\Delta t} \text{Id}_{N_{\mathrm{sp}}} & D_0^T & 0 & \dots & 0 \\
0 & \ddots & \ddots & \ddots & \vdots \\
\vdots & \dots & -\frac{1}{\Delta t} \text{Id}_{N_{\mathrm{sp}}} & D_{N_T-2}^T & 0 \\
0 & \dots & 0 & -\frac{1}{\Delta t} \text{Id}_{N_{\mathrm{sp}}} & D_{N_T-1}^T
\end{pmatrix}
\end{equation}

The identity block in the top-left enforces the initial condition $M^0 = m_0$, analogous to the bottom-right identity in $A_{\mathcal{U},\mathcal{U}}$, which enforces the terminal condition on $U$.
For non-separable Hamiltonians (Section~\ref{sec:numexp:nonlocalcongestion}), the sparsity pattern of $A_{\mathcal{M},\mathcal{M}}$ is unchanged but the off-diagonal blocks acquire additional contributions from $\partial_M \mathcal{B}$ beyond the transpose of $D_{n-1}$.
The remaining two Jacobian blocks are

\begin{equation} \label{eq:globalNewton-Jacobian_AUM_AMU}
A_{\mathcal{U}, \mathcal{M}} = \begin{pmatrix} 
0 & E_0 & 0 & \dots & 0 \\
0 & 0 & E_1 & \ddots & \vdots \\
\vdots & \ddots & \ddots & \ddots & 0 \\
0 & \dots & 0 & 0 & E_{N_T-1} \\
0 & \dots & 0 & 0 & E_{N_T}
\end{pmatrix}, \quad
A_{\mathcal{M}, \mathcal{U}} = \begin{pmatrix} 
0 & 0 & \dots & 0 & 0 \\
\tilde{E}_0 & 0 & \dots & 0 & 0 \\
0 & \tilde{E}_1 & \dots & 0 & 0 \\
\vdots & \ddots & \ddots & \vdots & \vdots \\
0 & \dots & 0 & \tilde{E}_{N_T-1} & 0
\end{pmatrix}
\end{equation}

Regarding $A_{\mathcal{M}, \mathcal{U}}$, all blocks in the first row are zero because the initial condition of $M$ does not depend on $U$. 
All blocks in the last column are zero because the terminal condition $U^{N_T}$ has no instantaneous effect on the forward flow of mass. 
Because the transport operator $\mathcal{B}_{i, j}(U^n, M^{n+1})$ in the discrete FP equation \eqref{eq:MFG-finitedifference-fb:FP} couples $U$ at time $n$ and $M$ at time $n+1$, $A_{\mathcal{M}, \mathcal{U}}$ has blocks on its first subdiagonal. 
The blocks $\tilde{E}_n = \nabla_U \mathcal{B}(U^n, M^{n+1})$ are highly sparse because of the structure of the transport operator \eqref{eq:transport-operator}. They represent the sensitivity of the mass distribution at time $n+1$ to changes in the optimal drift at time $n$.
In 1D, the operator uses a three-point stencil and therefore produces a tridiagonal matrix. 
In 2D, the operator uses the standard five-point stencil, producing a block-tridiagonal matrix with tridiagonal blocks, commonly described as pentadiagonal.

Regarding $A_{\mathcal{U},\mathcal{M}}$, the first column is all zeros since the backward value function is independent of the fixed initial mass.
The blocks $E_0, \ldots, E_{N_T - 1}$ lie on the superdiagonal because the HJB residual at time $n$ depends on $M^{n+1}$.
The last row is not all zeros: the terminal condition in $u$ depends on the terminal mass, placing the block $E_{N_T}$ on the main diagonal at position $(N_T, N_T)$.
For $n=0,\ldots,N_T-1$, any pointwise dependence of the Hamiltonian on the density makes $E_n$ diagonal, whether or not the Hamiltonian is additively separable in $m$ and $p$. The terminal block $E_{N_T}$ is diagonal when $g$ depends pointwise on the terminal density and vanishes when $g$ is independent of it. Nonlocal dependence in either the Hamiltonian or the terminal cost may instead produce dense blocks.
Algorithm~\ref{alg:global-newton} gives the pseudocode for the Newton method. 
The Newton method requires a good initial guess; see Section~\ref{sec:numericalmethods:practical} for the viscosity continuation strategy we use to obtain one.

The Picard method replaces the joint Newton solve with an outer Picard iteration over $(U, M)$. 
The iteration starts from an initial mass history $M^{(0)}$ built from the initial condition $m_0$ (held constant in time, or advanced by pure diffusion) and the value history $U^{(0)} \equiv 0$.
Each Picard iterate $k$ performs the following steps:

\textbf{(i) Backward HJB solve.} Given $M^{(k)}$, the nonlinear discrete HJB equation \eqref{eq:MFG-finitedifference-fb:HJB} is solved by marching backward from $n = N_T - 1$ to $n = 0$. At each timestep, Newton's method is initialized from the solution at the subsequent timestep $\hat{U}^{n+1}$; the residual and localized Jacobian are
\begin{align}
F_n(U) &= \frac{U - \hat{U}^{n+1}}{\Delta t} - \nu(\Delta_h U) + \tilde{H}(x, M^{(k),n+1}, [D_h U]), \\
J_n(U) &= \frac{1}{\Delta t} \mathrm{Id}_{N_{\mathrm{sp}}} - \nu \Delta_h + \sum_{l=1}^{2d} \mathrm{diag}\left(\frac{\partial \tilde{H}}{\partial q_l}(x, M^{(k),n+1}, [D_h U])\right) D_l,
\end{align}
where $J_n$ is tridiagonal in 1D and pentadiagonal in 2D.
The value function is then updated with a damping factor $\theta \in (0, 1]$:
\begin{equation*}
U^{(k+1)} = (1-\theta)U^{(k)} + \theta \hat{U}.
\end{equation*}

\textbf{(ii) Forward FP solve.} Given $U^{(k+1)}$, we solve the discrete FP equation~\eqref{eq:MFG-finitedifference-fb:FP} for a candidate $\hat{M}$ by marching forward from $M^0 = m_0$, one linear system per timestep. 
For locally separable Hamiltonians, $\partial \tilde{H} / \partial q_\ell$ is independent of $M$, so the transport operator $\mathcal{B}_{i,j}(U^{(k+1)}, \cdot)$ is already linear in $M^{n+1}$ and no further treatment is needed. 
For non-separable Hamiltonians, we evaluate $\partial \tilde{H} / \partial q_\ell$ at the previous Picard iterate $M^{(k)}$, restoring linearity in $M^{n+1}$ while introducing a lag handled by the outer Picard iteration. 
The mass is then updated with the same damping:
\begin{equation*}
M^{(k+1)} = (1-\theta)M^{(k)} + \theta \hat{M}.
\end{equation*}

The outer Picard iteration continues until the relative changes in both
$U$ and $M$ fall below a specified tolerance $\epsilon_P$.
Because $U^{(0)}\equiv0$, the relative convergence test is applied only
from the second Picard iteration onward.
Algorithm~\ref{alg:picard} gives the pseudocode for the Picard method.

We emphasize the substantial difference between the computational footprints of the two methods. 
Table~\ref{tab:jacobian-comparison} summarizes the Jacobian size, sparsity pattern, and number of linear solves per outer iteration for each method; the numerical experiments in Section~\ref{sec:numexp} supplement this structural comparison with wall-clock measurements at specific grid resolutions.

\begin{table}[h]
\centering
\begin{tabular}{lll}
\hline
 & Picard method & Newton method \\
\hline
Jacobian size & $N_{\mathrm{sp}}\times N_{\mathrm{sp}}$ & $2(N_T+1)N_{\mathrm{sp}}\times 2(N_T+1)N_{\mathrm{sp}}$ \\
Sparsity, 1D & tridiagonal & block-banded, tridiagonal sub-blocks \\
Sparsity, 2D (local Hamiltonian) & pentadiagonal & block-banded, pentadiagonal sub-blocks \\
Sparsity, 2D (nonlocal Hamiltonian) & pentadiagonal & dense sub-blocks in $A_{\mathcal{U},\mathcal{M}}$ and $A_{\mathcal{M},\mathcal{M}}$ \\
Solves per outer iteration & $\propto N_T \cdot (\text{Newton iters})$ & $1$ \\
\hline
\end{tabular}
\caption{Comparison of Jacobian structure and per-iteration cost for the Picard and Newton methods.}
\label{tab:jacobian-comparison}
\end{table}

The $N_{\mathrm{sp}}\times N_{\mathrm{sp}}$ systems arising in the Picard method are much smaller, but their efficiency depends on convergence of the outer Picard iteration. 
The Newton Jacobian is substantially larger, but Newton's method has local quadratic convergence when supplied with a suitable initial guess.
We have already noted advantages for the Picard method, such as the availability of a high-quality initial guess from temporal continuity and the use of damping to aid convergence. 
However, continuation in the viscosity parameter $\nu$ can improve the robustness of the Newton method. 
This technique helps overcome the limited basin of attraction for Newton's method and addresses the numerical difficulties inherent in the small-viscosity regime. Section~\ref{sec:numericalmethods:practical} examines these trade-offs and practical implementation details further.

Prior work helps explain why the Picard method has often been preferred in implementations.
Laurière~\cite{lauriere2021condcontrol} points out that solving the whole system by Newton's method can require costly full-matrix inversions in the presence of nonlocal terms, while Achdou and Laurière~\cite{achdou2020numericalaspects} emphasize the algebraic burden of differentiating both the Bellman and Kolmogorov equations.
These observations motivate the present comparison: modern symbolic and code-generation tools reduce the derivation barrier for the full Jacobian, leaving the computational trade-off between localized Picard solves and coupled Newton solves as the central numerical question.

\subsection{Practical considerations}
\label{sec:numericalmethods:practical}

\subsubsection{Continuation and initialization}

The Newton method requires an initial guess within the basin of attraction of the coupled space-time solve.
This is particularly restrictive at small viscosity, where the Laplacian provides less regularization and the Hamiltonian and transport terms dominate.
We therefore use continuation in the viscosity parameter: starting from a larger value $\nu_0$, where the naive initialization $U\equiv g_T$ and $M\equiv m_0$ is sufficient, we solve a sequence of problems with decreasing viscosities $\nu_0>\nu_1>\cdots>\nu_{\mathrm{target}}$, using each converged solution as the initial guess for the next solve.
The schedule is adaptive: attempted decreases are backtracked when the Newton solve fails.
The pseudocode is given in Algorithm~\ref{alg:newton_continuation}.

This continuation strategy is a numerical globalization heuristic.
Related continuation ideas are common in the numerical MFG literature~\cite{achdou2012planning,achdou2020numericalaspects}.

\subsubsection{Linear algebra backend}

We now discuss the specific linear-algebra backends used by these methods.
Concrete wall-clock measurements and grid-size thresholds are deferred to Section~\ref{sec:numexp}.
For the Picard method in one and two dimensions, the spatial Jacobians $J_n$ have size $N_{\mathrm{sp}}\times N_{\mathrm{sp}}$ and are tridiagonal and pentadiagonal, respectively.
In 1D, exact LU factorization introduces zero fill-in, allowing direct solvers---such as the \texttt{spsolve} routine from Python's \texttt{scipy.sparse.linalg}, which wraps the SuperLU library~\cite{Demmel1999superLU}---to remain computationally feasible for grids with $N_h$ on the order of $10^6$.
For pentadiagonal matrices, exact factorization does incur fill-in but remains feasible up to moderately large 2D grids.

In contrast, the Newton method requires solving a coupled space-time Jacobian.
For a two-dimensional grid with $200$ points in each spatial direction and $100$ timesteps, the system would have $2\times200^2\times101=8{,}080{,}000$ degrees of freedom, a scale at which direct factorization is intractable due to fill-in.
We deliberately sidestep this issue by keeping $N_T$ small enough that direct methods remain feasible: the two-dimensional experiments use $N_T=20$ and the hybrid timings $N_T=10$, so that every result reported here---for both methods---is obtained with the same direct sparse backend (\texttt{spsolve}/SuperLU) and the comparison is made under a common linear-algebra backend.
However, in $d=2$ the introduction of nonlocal congestion terms yields dense cross-coupling blocks within this space-time Jacobian, rendering explicit matrix assembly computationally intractable at the grid sizes of interest, which we elaborate on in Section~\ref{sec:numexp:nonlocalcongestion}.

\subsection{Hybrid method: \texorpdfstring{Picard-based acceleration of Newton continuation}{Picard-based acceleration of Newton continuation}}
\label{sec:numericalmethods:hybrid}

The hybrid method (Algorithm~\ref{alg:hybrid_strategy}) exploits the complementary cost-robustness profiles of the two solvers to reduce the number of expensive Newton solves needed to reach a target viscosity $\nu_{\text{target}}$ at a fixed target discretization $(N_h,N_T)$. We first use the Picard method on this target grid to descend in viscosity as far as convergence permits, lowering the damping parameter when needed. Once the Picard method fails, we hand off to the Newton method on the same grid and continue the remaining viscosity descent.

The design is motivated by the two-dimensional double-well sweep in Section~\ref{sec:numexp:double-well}: once the grid is sufficiently fine, further mesh refinement does not appear to lower the minimum viscosity reached by the Picard method in that benchmark. Thus the hybrid strategy tested here stays on the target grid throughout, rather than combining viscosity continuation with mesh refinement.

The strategy is organized in three phases:
\begin{enumerate}
    \item \textbf{Picard descent at the target discretization:}
    At the target discretization $(N_h, N_T)$, we apply the Picard method at decreasing viscosity levels starting from a large viscosity $\nu_{\text{start}}$. Every Picard attempt, including each retry with reduced damping, starts independently with $U^{(0)}\equiv0$ and a density history obtained by advancing $m_0$ under pure diffusion at the current viscosity. We denote this initialization by $\operatorname{InitPicard}(\nu)$ in Algorithm~\ref{alg:hybrid_strategy}. 
    In the experiments below, these viscosity levels are taken from the fixed schedule $\{1,0.3,0.1,0.03,0.01,0.003,0.001,0.0003,0.0001\}$. In Algorithm~\ref{alg:hybrid_strategy}, \texttt{NextViscosity} denotes the next smaller value in this list.
    We start with a rather large $\theta = 0.8$, and at each level, if the Picard method fails to converge, we decrease $\theta$ by $0.1$ and retry the same $\nu$ step, down to a floor of $\theta = 0.05$. 
    The damping parameter is retained across viscosity steps and never increased, reflecting the heuristic that lower viscosities require at least as much damping as higher ones. The converged Picard state is saved for a possible Newton handoff; it is not used to initialize the next Picard attempt.
    If the Picard method reaches $\nu_{\text{target}}$ unaided, the hybrid degenerates to the pure Picard method and Phases 2 and 3 are skipped.
    
    \item \textbf{Picard failure and identification of the handoff viscosity:}
    At the first viscosity where the Picard method fails after exhausting the damping schedule, we stop the Picard descent and retain the last converged solution.
    Its viscosity defines $\nu_{\text{limit}}$.
    Phase 2 records the failed Picard attempts at the next lower viscosity, which trigger the handoff to Newton.

    \item \textbf{Newton re-solve and continuation to $\nu_{\text{target}}$:}
    Starting from the last converged Picard solution, we first re-solve with Newton at $\nu_{\text{limit}}$ on the same discretization.
    This step is needed because the Picard method terminates based on the relative error between successive Picard iterates rather than the global residual; the re-solve reduces the full coupled residual before Newton continuation begins.
    Empirically, the converged Picard solution lies within the Newton basin of attraction at the same $(\nu, h)$.
    We then continue the viscosity descent to $\nu_{\text{target}}$ via the Newton method, using the previous converged $(U, M)$ as the initial guess at each step. 
\end{enumerate}

We record wall-clock times for each phase of the hybrid strategy using direct sparse linear solves, so that the phase timings identify the computational savings obtained by combining Picard and Newton.

\section{Numerical experiments}
\label{sec:numexp}

This section evaluates Picard-based acceleration of Newton continuation through a sequence of finite-difference MFG benchmarks, beginning with validation and an assessment of the strengths and limitations of the two component solvers.
The benchmark data are summarized in Table~\ref{tab:mfg_cases}. 
Unless otherwise stated, all computations use the semi-implicit finite-difference discretization described in Section~\ref{sec:numericalmethods:FD-discretization}; solver tolerances, damping parameters, continuation schedules, and failure criteria are reported in the corresponding experiment.

\begin{table}[h]
\centering
\caption{Summary of numerical benchmarks.}
\label{tab:mfg_cases}
\begin{tabular}{@{}llll@{}}
\toprule
\textbf{Case} & {\textbf{HJB nonlinearity}} & \textbf{Initial mass} $m_0(x)$ & \textbf{Terminal cost} $g(x)$ \\ \midrule
Discounted & {$\tfrac{1}{2}|p|^2 + \ln m$} & {$\bar m(x)$} & {$\bar u(x)$} \\
Ergodic & {$\tfrac{1}{2}|p|^2 + \ln m$} & $m_0 \equiv \bar m(\pm L)$ & {$\bar u(x)$} \\
Congestion (nonlocal) & $\frac{|p|^2}{2(1+4\bar m)^\beta} - \zeta \bar m$ & Square wave & Double-quadratic well \\ 
Traffic light & $\frac{1}{2}|p|^2 - (V(x, t) + \kappa m^{\alpha})$ & Gaussian (centered at left) &  $5(x-1)^2$\\
Double-well & $\frac{1}{2}|p|^2 - (V(x) + \kappa m^\alpha)$ & Gaussian at origin & $V(x)$---quartic potential \\
\bottomrule
\end{tabular}
\end{table}

The first two examples are validation tests. 
The discounted problem compares both solvers against a known stationary solution on a truncated domain. 
The ergodic example uses a manufactured stationary structure to study the vanishing-discount limit. 
The next two examples identify regimes favoring one method over the other: nonlocal congestion favors Picard because of the cost of assembling and solving the coupled Newton systems, while the traffic-light problem favors Newton continuation because of its greater convergence robustness in difficult low-viscosity regimes with temporal shocks. The double-well problem then illustrates the hybrid strategy using direct linear solves.

\subsection{A discounted problem} 
\label{sec:numexp:discounted}

Our first validation example is the stationary discounted problem \eqref{eq:stationary-discounted-PDE-system}, for which Gueant \cite{Gueant2009reference} derived a closed-form solution on $\mathbb{R}^d$. 
A family of solutions $(\bar{u}, \bar{m})$, parameterized by its mean $\mu \in \mathbb{R}^d$, is given by
\begin{equation} \label{eq:discounted-analytical-u-m}
\bar{u}(x) = \eta |x - \mu|^2 + \omega, \quad \bar{m}(x) = \mathcal{N}(\mu, s^2) =  \frac{1}{(2\pi s^2)^{d/2}} \exp\left(-\frac{|x - \mu|^2}{2s^2}\right),
\end{equation}
and the constants $s^2 > 0$, $\eta > 0$, and $\omega$ are
\begin{equation} \label{eq:discounted-analytical-constants}
s^2 = \frac{\sigma^4}{4 - 2\rho\sigma^2}, \quad \eta = \frac{\sigma^2}{4s^2} = \frac{1}{\sigma^2} - \frac{\rho}{2}, \quad \omega = \frac{1}{\rho} \left[ \eta d \sigma^2 - \frac{d}{2} \ln \left( \frac{2\eta}{\pi \sigma^2} \right) \right].
\end{equation}
These formulas give the analytic stationary reference for this validation test. They require $\rho\sigma^2<2$, which guarantees $s^2>0$ and $\eta>0$; all parameter choices below satisfy this condition. 

To adapt this infinite-domain solution to a bounded computational grid, we truncate the domain and enforce Dirichlet boundary conditions derived from the values of the exact stationary solution.
We center the domain at a chosen $\mu$ and impose fixed Dirichlet data---obtained directly from the closed-form formulas---for all $t \in [0,T]$.
We initialize $m_0=\bar m$ and $u_T=\bar u$ and impose Dirichlet boundary conditions consistent with the same stationary solution.
This makes the experiment a direct reference-validation test for the finite-difference discretization and nonlinear solvers.
For simplicity, we take $\mu=0$ and $\Omega=[-L,L]$.

\begin{figure}[h]
    \centering
\begin{tikzpicture}
\begin{axis}[
    scale only axis,
    width=6.75cm,
    height=4.5cm,
    title={Density Comparison ($T=10$)},
    xlabel={$x$},
    ylabel={$m(x)$},
    ytick={0, 0.2, 0.4, 0.6, 0.8},  
    grid=major,
    grid style={dashed, gray!60},
    legend pos=north east,
    legend cell align={left},
legend style={
        at={(0.98,0.98)},      
        anchor=north east,    
        font=\scriptsize,     
        nodes={scale=0.8, transform shape}, 
        row sep=0.1pt,
        legend cell align={left}
    },
    tick label style={font=\small},
    label style={font=\small},
    title style={font=\scriptsize\bfseries}
]

\addplot[color=black, dashed, thick] table[x index=0, y index=1] {tikz/solution_slices/discounted/discounted-1d-picard-density.dat};
\addlegendentry{Stationary $\overline{m}$}

\addplot[color=blue, thick] table[x index=0, y index=2] {tikz/solution_slices/discounted/discounted-1d-picard-density.dat};
\addlegendentry{Numerical $m(T/2)$}

\end{axis}
\end{tikzpicture}    
\begin{tikzpicture}
\begin{axis}[
    scale only axis,
    width=6.75cm,
    height=4.5cm,
    title={Value Comparison ($T=10$)},
    xlabel={$x$},
    ylabel={$u(x)$},
    ytick={2, 3, 4, 5, 6, 7},    
    grid=major,
    grid style={dashed, gray!60},
    legend pos=north east,
    legend cell align={left},
legend style={
        at={(0.7,0.98)},      
        anchor=north east,    
        font=\scriptsize,     
        nodes={scale=0.8, transform shape}, 
        row sep=0.1pt,
        legend cell align={left}
    },
    tick label style={font=\small},
    label style={font=\small},
    title style={font=\scriptsize\bfseries}
]

\addplot[color=black, dashed, thick] table[x index=0, y index=1] {tikz/solution_slices/discounted/discounted-1d-picard-value-function.dat};
\addlegendentry{Stationary $\overline{u}$}

\addplot[color=red, thick] table[x index=0, y index=2] {tikz/solution_slices/discounted/discounted-1d-picard-value-function.dat};
\addlegendentry{Numerical $u(T/2)$}

\end{axis}
\end{tikzpicture}
 \caption{The discounted problem in $d = 1$. 
 Numerical solution of the forward-backward problem at time $T/2$ 
 compared to the closed-form solution for the stationary problem for $m$ \textit{(Left)} and $u$ \textit{(Right)}.} 
 \label{fig:discounted-picard-solution-slices}
\end{figure}
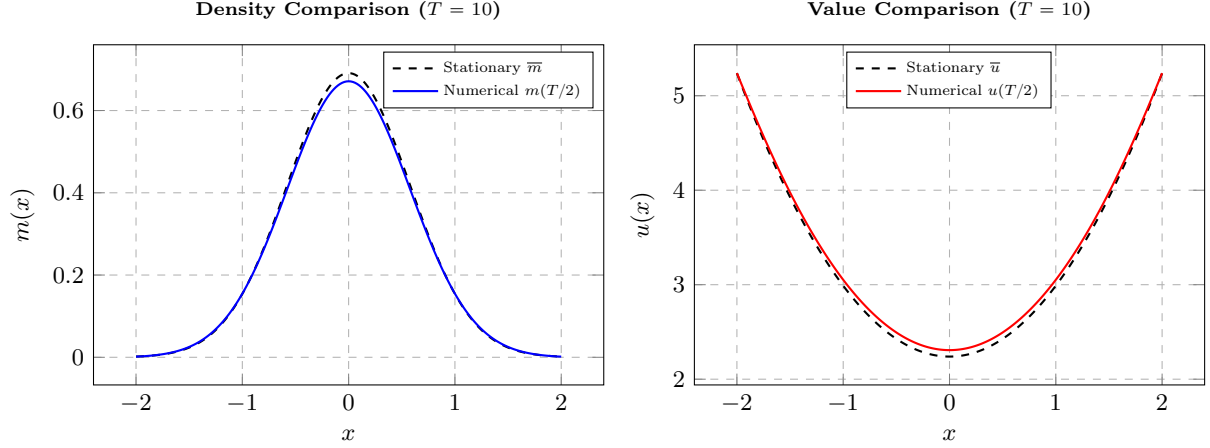

{The variance $s^2=\nu^2/(1-\rho\nu)$ of $\bar{m}$ scales as $\nu^2$ for fixed $\rho$.} 
As $\nu$ decreases with $L$ fixed, the Gaussian density $m$ approaches zero outside an increasingly narrow neighborhood of the origin. 
Resolving this sharp concentration requires additional grid refinement near the center $\mu$.
To ensure a well-conditioned system while capturing the essential physics, we select the parameters $\rho = 0.5$ and $\nu = 0.5$ on a domain with spatial half-width $L = 2$. 
We retain the time horizon $T = 10$, and the space-time discretization is set to $N_h = 200$ and $N_T = 200$.
{Figure~\ref{fig:discounted-picard-solution-slices} compares the Picard solution of the corresponding finite-horizon problem at time $T/2$ with the analytic stationary reference \eqref{eq:discounted-analytical-u-m}.}

\begin{figure}[H]
    \centering
\begin{tikzpicture}
\begin{semilogyaxis}[
    scale only axis,
    width=6.75cm,
    height=4.5cm,
    title={Picard Iteration History with $\theta = 0.8$},
    xlabel={Picard Iterate},
    ylabel={Relative Error (Log Scale)},
    xmin=0, xmax=31,
    ymin=1e-7, ymax=2,
    xtick={5,10,15,20,25,30},
    ytick={1e-6, 1e-5, 1e-4, 1e-3, 1e-2, 1e-1, 1},
    grid=both,
    grid style={line width=.1pt, draw=gray!10},
    major grid style={line width=.2pt, draw=gray!50},
    legend pos=north east,
    tick label style={font=\small},
    label style={font=\small},
    title style={font=\scriptsize\bfseries}
]

\addplot[
    color=red,
    mark=*,
    mark size=1.5pt,
    thick
] coordinates {
    (1,1.0000e+00) (2,1.3623e-01) (3,7.9176e-02) (4,1.1506e-01) (5,9.1495e-02)
    (6,5.1919e-02) (7,2.4578e-02) (8,1.0694e-02) (9,4.5954e-03) (10,2.1112e-03)
    (11,1.1127e-03) (12,6.6725e-04) (13,4.2708e-04) (14,2.7884e-04) (15,1.8247e-04)
    (16,1.1913e-04) (17,7.7549e-05) (18,5.0369e-05) (19,3.2663e-05) (20,2.1159e-05)
    (21,1.3697e-05) (22,8.8627e-06) (23,5.7330e-06) (24,3.7078e-06) (25,2.3977e-06)
    (26,1.5504e-06) (27,1.0025e-06) (28,6.4815e-07) (29,4.1905e-07) (30,2.7092e-07)
};
\addlegendentry{Rel Error $u$}

\addplot[
    color=blue,
    mark=square*,
    mark size=1.5pt,
    thick
] coordinates {
    (1,3.3963e-01) (2,2.9916e-01) (3,2.9396e-01) (4,2.2460e-01) (5,1.2244e-01)
    (6,6.0550e-02) (7,3.0787e-02) (8,1.6680e-02) (9,9.6077e-03) (10,5.7894e-03)
    (11,3.5905e-03) (12,2.2653e-03) (13,1.4436e-03) (14,9.2535e-04) (15,5.9514e-04)
    (16,3.8353e-04) (17,2.4744e-04) (18,1.5976e-04) (19,1.0319e-04) (20,6.6663e-05)
    (21,4.3074e-05) (22,2.7834e-05) (23,1.7987e-05) (24,1.1624e-05) (25,7.5119e-06)
    (26,4.8546e-06) (27,3.1373e-06) (28,2.0274e-06) (29,1.3102e-06) (30,8.4672e-07)
};
\addlegendentry{Rel Error $m$}

\draw[dashed, thick] (axis cs:0, 1e-6) -- (axis cs:31, 1e-6) node[pos=0.05, above right] {$\epsilon_P = 10^{-6}$};

\end{semilogyaxis}
\end{tikzpicture}
    \hfill
\begin{tikzpicture}
\begin{axis}[
    scale only axis,
    width=6.75cm,
    height=4.5cm,
    title={Newton Iterations per Picard Iteration},
    xlabel={Picard Iterate},
    ylabel={Number of Newton Iterations},
    xmin=0, xmax=31,
    ymin=550, ymax=650,
    xtick={5,10,15,20,25,30},
    grid=both,
    grid style={line width=.1pt, draw=gray!10},
    major grid style={line width=.2pt, draw=gray!50},
    legend pos=north east,
    tick label style={font=\small},
    label style={font=\small},
    title style={font=\scriptsize\bfseries}
]
\addplot[
    color=blue,
    mark=square*,
    mark size=1.5pt,
    thick
] coordinates {
    (1,635) (2,631) (3,630) (4,619) (5,600)
    (6,600) (7,600) (8,596) (9,582) (10,563)
    (11,579) (12,585) (13,588) (14,589) (15,589)
    (16,590) (17,590) (18,591) (19,591) (20,591)
    (21,591) (22,591) (23,591) (24,591) (25,591)
    (26,591) (27,591) (28,591) (29,591) (30,591)
};
\addlegendentry{Newton Its}
\end{axis}
\end{tikzpicture}
    \caption{\textit{(Left)} Picard iteration history for the Picard method applied to the discounted problem with $\theta = 0.8$. \textit{(Right)} Total number of Newton iterations in the backward HJB solve per Picard iterate.} 
    \label{fig:discounted-picard-history}
\end{figure}
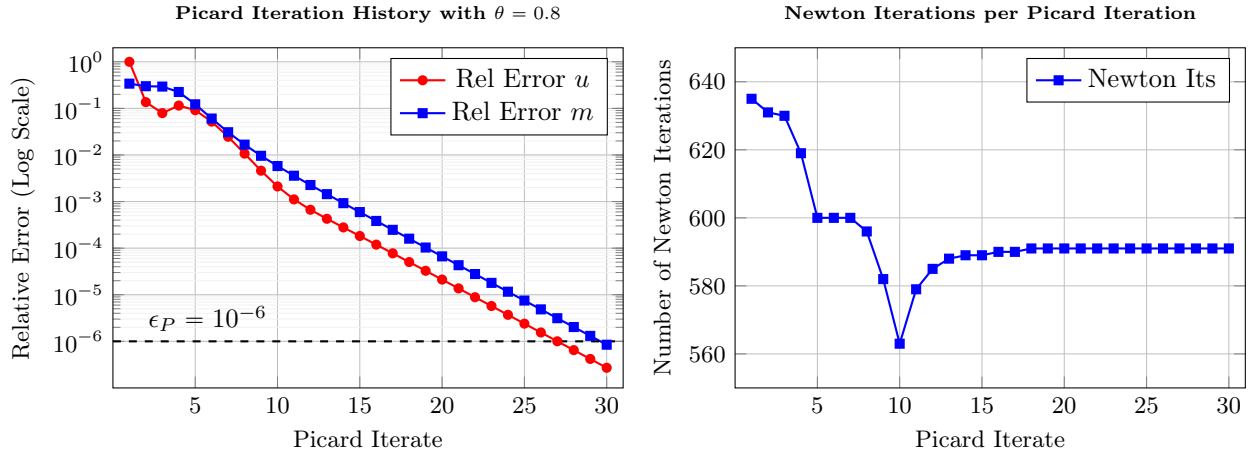

The left panel of Figure~\ref{fig:discounted-picard-history} shows a typical well-behaved Picard iteration. 
{After a short transient, the relative error for successive iterations of $u$ and $m$ decreases log-linearly and reaches $\epsilon_P$ in 30 iterations.}
As a rule of thumb, monotonically decreasing relative errors indicate that the damping parameter $\theta$ may be increased.
{Within each Picard iteration, the backward HJB sweep requires an average of approximately 2.8--3.2 Newton iterations per timestep to reach the tolerance $\epsilon_N$.} This low iteration count is expected because the value computed at the subsequent timestep provides a high-quality initial guess.
{Aggregated over $N_T = 200$ timesteps, this gives 563--635 Newton iterations per Picard iterate (Figure~\ref{fig:discounted-picard-history}, right).}

\begin{figure}[H]
    \centering
\begin{tikzpicture}
\begin{axis}[
    scale only axis,
    width=6.75cm,
    height=4.5cm,
    ymin=0.4,
    ymax=1.1,
    ytick={0.4, 0.6, 0.8, 1.0},
    xmin=0,
    xmax=36,
    xtick={5,10,15,20,25,30,35},
    xlabel={Continuation Attempt},
    ylabel={$\nu$ (Viscosity)},
    title={Viscosity Continuation for Global Newton},
    grid=major,
    grid style={dashed, gray!60},
    legend pos=north east,
    legend cell align={left},
    legend style={font=\small, nodes={scale=0.9, transform shape}},
    tick label style={font=\small},
    label style={font=\small},
    title style={font=\scriptsize\bfseries}
]
\addplot[
    only marks,
    mark=*,
    mark size=2pt,
    color=black,
    fill=black
] coordinates {
    (1,1.0000) (4,0.8750) (8,0.8281) (12,0.7871) (16,0.7512)
    (21,0.7355) (25,0.7061) (29,0.6803) (32,0.6352) (34,0.5676)
    (35,0.5000)
};
\addlegendentry{Converged}

\addplot[
    only marks,
    mark=x,
    mark size=3pt,
    thick,
    color=red
] coordinates {
    (2,0.5000) (3,0.7500) (5,0.5000) (6,0.6875) (7,0.7812)
    (9,0.5000) (10,0.6641) (11,0.7461) (13,0.5000) (14,0.6436)
    (15,0.7153) (17,0.5000) (18,0.6256) (19,0.6884) (20,0.7198)
    (22,0.5000) (23,0.6178) (24,0.6766) (26,0.5000) (27,0.6030)
    (28,0.6546) (30,0.5000) (31,0.5902) (33,0.5000)
};
\addlegendentry{Failed}

\addplot[
    color=black,
    dashed,
    domain=0:36,
    samples=2,
    thick
] {0.5};
\addlegendentry{Target $\nu = 0.5$}

\addplot[
    color=gray,
    draw opacity=0.5,
    thick,
    forget plot
] coordinates {
    (1,1.0000) (2,0.5000) (3,0.7500) (4,0.8750) (5,0.5000)
    (6,0.6875) (7,0.7812) (8,0.8281) (9,0.5000) (10,0.6641)
    (11,0.7461) (12,0.7871) (13,0.5000) (14,0.6436) (15,0.7153)
    (16,0.7512) (17,0.5000) (18,0.6256) (19,0.6884) (20,0.7198)
    (21,0.7355) (22,0.5000) (23,0.6178) (24,0.6766) (25,0.7061)
    (26,0.5000) (27,0.6030) (28,0.6546) (29,0.6803) (30,0.5000)
    (31,0.5902) (32,0.6352) (33,0.5000) (34,0.5676) (35,0.5000)
};
\end{axis}
\end{tikzpicture}
    \hfill
\begin{tikzpicture}
\begin{axis}[
    scale only axis,
    width=6.75cm,
    height=4.5cm,
    title={Newton Iterations per Continuation Attempt},
    xlabel={Continuation Attempt},
    ylabel={Number of Newton Iterations},
    xmin=0, xmax=36,
    ymin=0, ymax=28,
    xtick={5,10,15,20,25,30,35},
    ytick={0,5,10,15,20},
    grid=major,
    grid style={dashed, gray!60},
    legend pos=north east,
    legend cell align={left},
    legend style={font=\small, nodes={scale=0.9, transform shape}},
    tick label style={font=\small},
    label style={font=\small},
    title style={font=\scriptsize\bfseries}
]
\addplot[
    color=gray,
    draw opacity=0.5,
    thick,
    forget plot
] coordinates {
    (1,6) (2,20) (3,20) (4,5) (5,20)
    (6,20) (7,20) (8,4) (9,20) (10,20)
    (11,20) (12,4) (13,20) (14,20) (15,20)
    (16,6) (17,20) (18,20) (19,20) (20,20)
    (21,4) (22,20) (23,20) (24,20) (25,5)
    (26,20) (27,20) (28,20) (29,4) (30,20)
    (31,20) (32,5) (33,20) (34,5) (35,5)
};

\addplot[
    only marks,
    mark=*,
    mark size=2pt,
    color=black,
    fill=black
] coordinates {
    (1,6) (4,5) (8,4) (12,4) (16,6)
    (21,4) (25,5) (29,4) (32,5) (34,5)
    (35,5)
};
\addlegendentry{Converged}

\addplot[
    only marks,
    mark=x,
    mark size=3pt,
    thick,
    color=red
] coordinates {
    (2,20) (3,20) (5,20) (6,20) (7,20)
    (9,20) (10,20) (11,20) (13,20) (14,20)
    (15,20) (17,20) (18,20) (19,20) (20,20)
    (22,20) (23,20) (24,20) (26,20) (27,20)
    (28,20) (30,20) (31,20) (33,20)
};
\addlegendentry{Failed}
\end{axis}
\end{tikzpicture}
    \caption{Newton continuation history for the discounted problem. \textit{(Left)} {The 35 viscosity-continuation attempts, of which 24 fail, including nine direct attempts at $\nu=0.5$.} \textit{(Right)} Number of Newton updates, equivalently global linear solves, performed in each attempt. {The final attempt reaches $\nu=0.5$.}} 
    \label{fig:discounted-newton-history}
\end{figure}
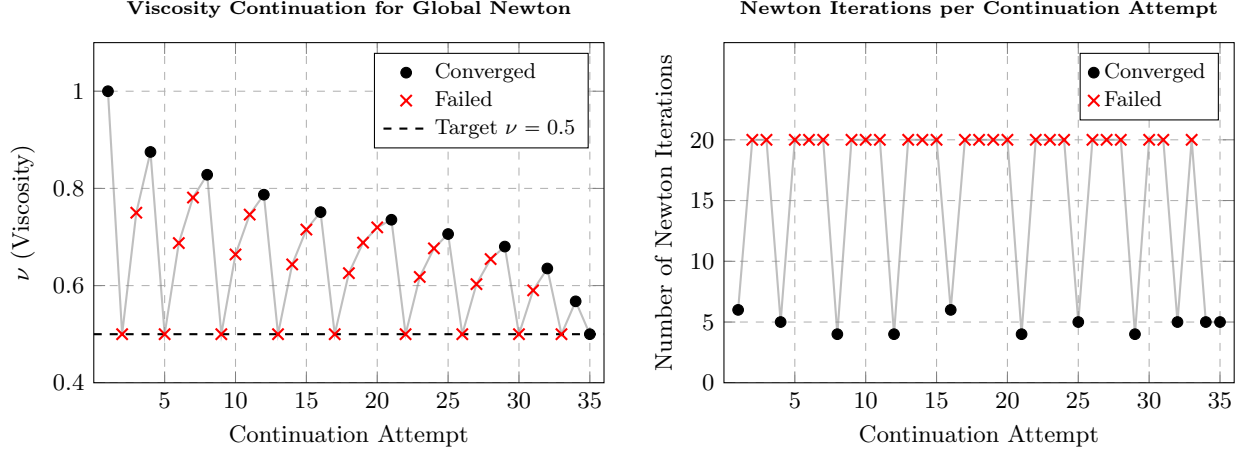

The left panel of Figure~\ref{fig:discounted-newton-history} shows the continuation path from $\nu=1$ to the target value $\nu=0.5$, using a geometric reduction factor $\tau=0.5$ and backtracking factor $\gamma=0.5$ (Algorithm~\ref{alg:newton_continuation}).
Every direct attempt to step to $\nu=0.5$ diverges until the current viscosity is close enough, and the backtracking produces 35 continuation attempts in total, of which 11 converge. 
Each converged value of $\nu$ requires 4--6 Newton updates.
The backtracking schedule is certainly not optimal, but chosen here to show its necessity. 

\begin{figure}[h]
    \centering
\newcommand{\tridiagbox}[3]{
    \draw[gray!40, thin] (#1, #2) rectangle (#1+1, #2+1);
    \draw[#3, thick] (#1+0.2, #2+0.8) -- (#1+0.8, #2+0.2);
    \draw[#3, thick] (#1+0.4, #2+0.8) -- (#1+0.8, #2+0.4);
    \draw[#3, thick] (#1+0.2, #2+0.6) -- (#1+0.6, #2+0.2);
}
\newcommand{\diagbox}[3]{
    \draw[gray!40, thin] (#1, #2) rectangle (#1+1, #2+1);
    \draw[#3, thick] (#1+0.2, #2+0.8) -- (#1+0.8, #2+0.2);
}

\begin{tikzpicture}[scale=0.6]
    \begin{scope}[shift={(0,3)}] 
        \draw[blue, thick] (0.4,3.6) -- (3.6,0.4); 
        \draw[blue, thick] (1.0,3.6) -- (3.6,1.0); 
        \draw[blue, thick] (0.4,3.0) -- (3.0,0.4); 
        \draw[black, thick] (0,0) rectangle (4,4);
        
        \node at (2, 4.6) {$\mathbf{J_n(U)}$};
        \node at (2, -0.8) {\textbf{Picard Method}};
        \node at (2, -1.6) {\textbf{(1D Separable Hamiltonian)}};
    \end{scope}

    \begin{scope}[shift={(8,0)}] 
        \foreach \i in {0,1,2,3,4} { \tridiagbox{\i}{9-\i}{blue} }
        \foreach \i in {1,2,3,4} { \diagbox{\i}{10-\i}{blue!40} }

        \foreach \i in {0,1,2,3,4} { \tridiagbox{5+\i}{4-\i}{blue} }
        \foreach \i in {0,1,2,3} { \diagbox{5+\i}{3-\i}{blue!40} }

        \foreach \i in {0,1,2,3} { \tridiagbox{\i}{3-\i}{red} }

        \foreach \i in {1,2,3,4} { \diagbox{5+\i}{10-\i}{green!60!black} }
        \diagbox{9}{5}{green!60!black}

        \draw[dashed, thick] (5,0) -- (5,10);
        \draw[dashed, thick] (0,5) -- (10,5);
        \draw[black, thick] (0,0) rectangle (10,10);
        
        \node at (2.5, 10.6) {$\mathbf{A_{\mathcal{U},\mathcal{U}}}$};
        \node at (7.5, 10.6) {$\mathbf{A_{\mathcal{U},\mathcal{M}}}$};
        \node at (2.5, -0.6) {$\mathbf{A_{\mathcal{M},\mathcal{U}}}$};
        \node at (7.5, -0.6) {$\mathbf{A_{\mathcal{M},\mathcal{M}}}$};
        \node at (5, -1.3) {\textbf{Global Newton Method (1D Separable Hamiltonian)}};
    \end{scope}
\end{tikzpicture}
    \caption{Sparsity structure for the 1D separable Hamiltonian. \textit{(Left)} The Picard-method Jacobian $J_n(U)$, which is strictly tridiagonal of size $(N_h+1) \times (N_h+1)$. \textit{(Right)} The Newton-method space-time Jacobian of size $2(N_T+1)(N_h+1) \times 2(N_T+1)(N_h+1)$. 
    Each block represents an $(N_h+1) \times (N_h+1)$ spatial operator. 
    Note that $D_i$ and $\tilde{E}_i$ are strictly tridiagonal, while $E_i$ remains diagonal due to local separability.}
\label{fig:sparsity_separableH_1d}
\end{figure}
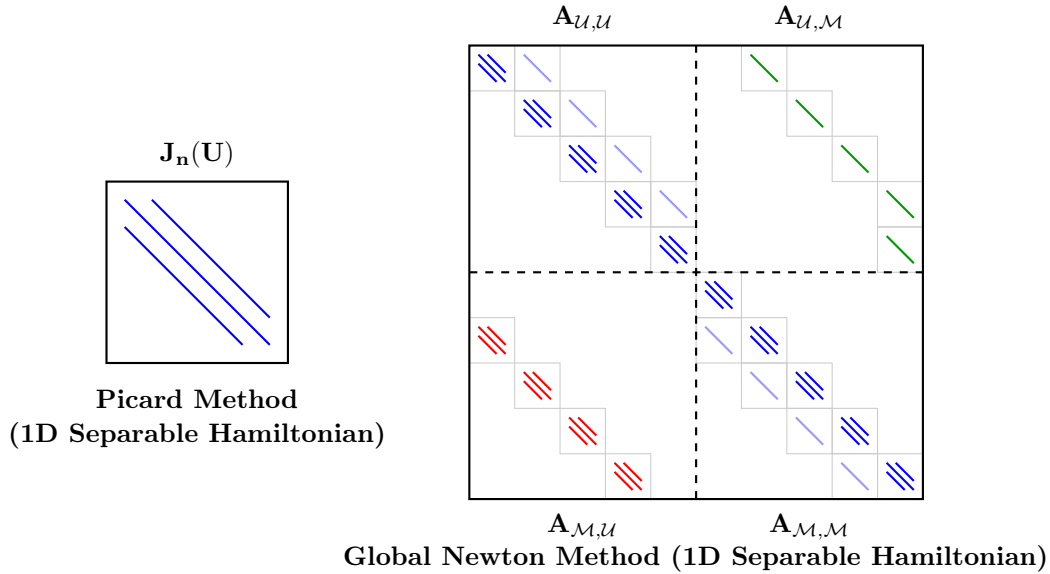

We next compare the linear-solve costs of the two methods in $d = 1$. 
According to the sparsity patterns in Figure~\ref{fig:sparsity_separableH_1d}, the Picard method requires repeated solves of a tridiagonal $(N_h+1) \times (N_h+1)$ matrix, which for $N_h = 200$ spatial intervals is a $201 \times 201$ system. 
The Newton method assembles the four block matrices introduced in Section~\ref{sec:numericalmethods:Newton} into a space-time Jacobian of size $2(N_T+1)(N_h+1) \times 2(N_T+1)(N_h+1)$, which for $N_h = N_T = 200$ is $80{,}802 \times 80{,}802$.
The blocks of $D_i$ in $A_{\cU, \cU}$ and $A_{\cM, \cM}$ are tridiagonal, as are the blocks of $\tilde{E}_i$ in $A_{\cM, \cU}$. 
The blocks in $A_{\cU, \cM}$ are diagonal due to local separability. 

The disparity in computational cost between the two methods is driven by the structural properties of their respective linear systems. 
The total number of spatial solves includes those from both the backward HJB sweep and the forward Fokker--Planck (FP) sweep. 
{According to the convergence history, the HJB solve averages roughly 2.97 inner Newton iterations per timestep, plus one additional solve required to march the linear FP equation forward. The Picard loop converges in 30 iterations (Figure~\ref{fig:discounted-picard-history}, left).
The total number of solves in the experiment is
$$(2.9748 \text{ HJB Newton iters} + 1 \text{ FP step}) \times N_T (200) \times \text{Picard iterations (30)} = 23{,}849 \text{ solves}.$$
The total execution time for the entire Picard process is 14.60 seconds.}
{In contrast, the Newton method with continuation requires 35 continuation attempts (24 failures triggering backtracks) and 533 global Newton linear solves, of which 53 occur in the 11 converged attempts, for a wall-clock total of approximately 2007.8 seconds (253.3 seconds in the converged attempts).} 
The backtracking strategy in $\nu$ is deliberately simple and not optimized; nevertheless, the timing comparison shows the cost gap between the coupled Newton solve and the Picard method in this benchmark.

{Thus, in this one-dimensional benchmark, the sequential temporal marching and reduced-order spatial Jacobians of the Picard method yield a roughly 137-fold speedup (17-fold if the failed backtracks are excluded) over the coupled space-time Newton solve.}
Section~\ref{sec:numexp:ergodic} next stress-tests both methods in the vanishing-discount limit $\rho \to 0$, where the additive level of the value function becomes unbounded and Newton continuation must traverse several orders of magnitude in $\rho$.

\subsection{An ergodic problem} 
\label{sec:numexp:ergodic}

{This subsection contains two related validation checks.
The first is an unforced stationary ergodic reference used to check the finite-horizon Picard solver against a known profile.}
The second is a discounted $\rho$-continuation experiment used to test the discrete vanishing-discount relation $\rho u_\rho \to \lambda_{\Omega,h}$ over four orders of magnitude in $\rho$.
The spatial operators are the same as in the discounted benchmark, so we do not repeat the per-iteration sparsity and cost discussion.

{For the stationary reference check, the three constants in \eqref{eq:discounted-analytical-constants} have the following limits as $\rho\to0$:}

\begin{equation} 
\label{eq:discounted-analytical-constants-rholimit}
\bar s^2 = \frac{\sigma^4}{4}, \quad \bar \eta = \frac{1}{\sigma^2}, \quad \omega \rightarrow \infty.
\end{equation}
{Substituting the limiting Gaussian--quadratic profiles into the stationary ergodic system \eqref{eq:stationary-ergodic-PDE-system}, the product $\rho u_\rho(x)$ approaches the ergodic constant}
\begin{equation} \label{eq:true-ergodic-constant}
\lambda= \bar \eta d \sigma^2 - \frac{d}{2} \ln \left( \frac{2\bar \eta}{\pi \sigma^2} \right)  = d - \frac{d}{2}\ln\left(\frac{2}{\pi\sigma^4}\right).
\end{equation}

The divergent $\omega$ reflects the $\mathcal{O}(1/\rho)$ blow-up of the discounted value function as $\rho \to 0$; we discard it and define a new constant $\bar \omega$, fixed below by the ergodic uniqueness condition $\int_\Omega \bar u(x)\, dx = 0$.
For consistency with Section~\ref{sec:numexp:discounted}, we pick $\mu$ as the origin.
We choose the following form for the stationary mass:
\begin{equation} \label{eq:ergodic-analytical-m}
\bar m(x) = \frac{1}{(2\pi \bar s^2)^{d/2}} \exp\left(-\frac{|x - \mu|^2}{2\bar s^2}\right).
\end{equation}
The corresponding value function satisfying the FP equation is
\begin{equation} \label{eq:ergodic-analytical-u}
    \bar u(x) = \bar \eta |x - \mu|^2 + \bar \omega.
\end{equation}
{On the computational grid, we normalize the sampled density so that $\sum_i \bar m_i\Delta x=1$ and subtract the discrete spatial mean from $\bar u$. The profiles \eqref{eq:ergodic-analytical-m}--\eqref{eq:ergodic-analytical-u}, together with \eqref{eq:true-ergodic-constant}, provide the analytic reference for this unforced problem.}

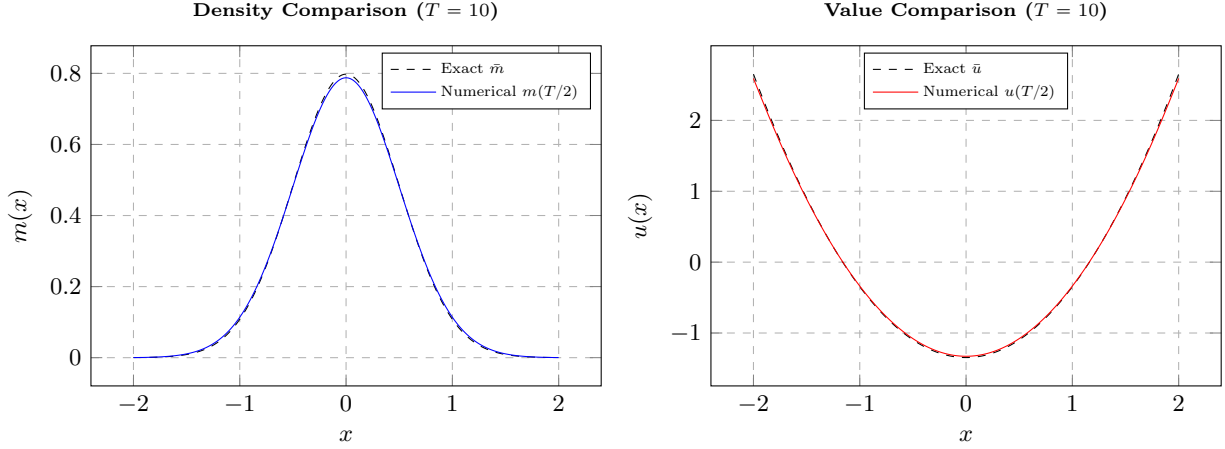
\begin{figure}[H]
    \centering
\begin{tikzpicture}
\begin{axis}[
    scale only axis,
    width=6.75cm,
    height=4.5cm,
    title={Density Comparison ($T=10$)},
    xlabel={$x$},
    ylabel={$m(x)$},
    ytick={0, 0.2, 0.4, 0.6, 0.8},
    grid=major,
    grid style={dashed, gray!60},
    legend pos=north east,
    legend cell align={left},
legend style={
        at={(0.98,0.98)},      
        anchor=north east,    
        font=\scriptsize,     
        nodes={scale=0.8, transform shape}, 
        row sep=0.1pt,
        legend cell align={left}
    },
    tick label style={font=\small},
    label style={font=\small},
    title style={font=\scriptsize\bfseries}
]
\addplot[black, dashed] table[x index=0, y index=1] {tikz/solution_slices/ergodic/ergodic-1d-picard-density.dat};
\addlegendentry{Exact $\bar{m}$}

\addplot[blue] table[x index=0, y index=2] {tikz/solution_slices/ergodic/ergodic-1d-picard-density.dat};
\addlegendentry{Numerical $m(T/2)$}
\end{axis}
\end{tikzpicture}    
\begin{tikzpicture}
\begin{axis}[
    scale only axis,
    width=6.75cm,
    height=4.5cm,
    title={Value Comparison ($T=10$)},
    ytick={-1, 0, 1, 2},
    xlabel={$x$},
    ylabel={$u(x)$},
    grid=major,
    grid style={dashed, gray!60},
    legend pos=north east,
    legend cell align={left},
legend style={
        at={(0.7,0.98)},      
        anchor=north east,    
        font=\scriptsize,     
        nodes={scale=0.8, transform shape}, 
        row sep=0.1pt,
        legend cell align={left}
    },
    tick label style={font=\small},
    label style={font=\small},
    title style={font=\scriptsize\bfseries}
]
\addplot[black, dashed] table[x index=0, y index=1] {tikz/solution_slices/ergodic/ergodic-1d-picard-value-function.dat};
\addlegendentry{Exact $\bar{u}$}

\addplot[red] table[x index=0, y index=2] {tikz/solution_slices/ergodic/ergodic-1d-picard-value-function.dat};
\addlegendentry{Numerical $u(T/2)$}
\end{axis}
\end{tikzpicture}
 \caption{The ergodic problem in $d = 1$. 
 Numerical solution of the forward-backward problem at time $T/2$ 
 compared to the closed-form solution for the stationary problem for $m$ \textit{(Left)} and $u$ \textit{(Right)}.} 
 \label{fig:ergodic-picard-solution-slices}
\end{figure}

{Using the same spatial and temporal parameters as in Section~\ref{sec:numexp:discounted}, we solve the corresponding finite-horizon problem, with the ergodic constant $\lambda$ from \eqref{eq:true-ergodic-constant} inserted in the HJB equation, Dirichlet boundary values taken from $(\bar u,\bar m)$, a flat initial density $m_0\equiv\bar m(\pm L)$ (renormalized to unit mass), and terminal value $u_T=\bar u$, and extract the numerical solution at $T/2$.
We anchor the terminal value to $\bar u$ rather than to a constant: because the coupling is not monotone, the finite-horizon problem started from a constant terminal value converges to a different forward-backward solution whose $T/2$ slice is not the stationary state, and this persists when the horizon is extended to $T=40$.}
Figure~\ref{fig:ergodic-picard-solution-slices} compares the numerical solution with the closed-form stationary solution.

{In contrast to the discounted case, the ergodic value function is defined only up to an additive constant. For the plotted comparison, we subtract the discrete spatial mean from each value profile.}

We now turn from the manufactured stationary check to the vanishing-discount experiment announced at the start of this section.
The continuous relation discussed in Section~\ref{sec:eqnsanddisc:connections} is asymptotic and stationary; on the truncated computational grid, we instead test its discrete analogue.
The target is therefore not the continuum constant in \eqref{eq:true-ergodic-constant}, but the discrete constant $\lambda_{\Omega,h}$ computed on the same domain and grid.

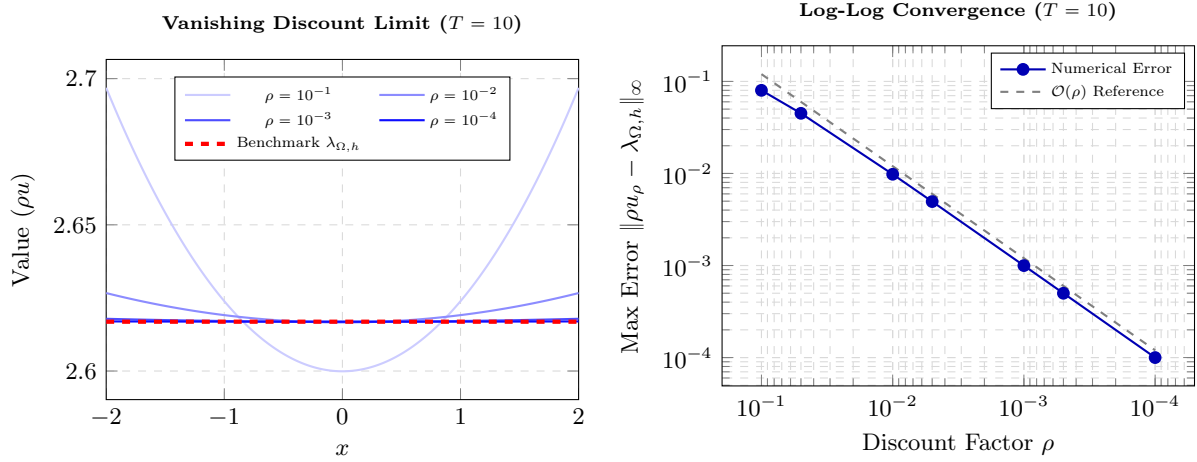
\begin{figure}[H]
    \centering
\begin{tikzpicture}
\begin{axis}[
    scale only axis,
    width=6.25cm,
    height=4.5cm,
    title={Vanishing Discount Limit ($T=10$)},
    xlabel={$x$},
    ylabel={Value ($\rho u$)},
    grid=major,
    grid style={dashed, gray!30},
    legend style={
        at={(0.5,0.95)},
        anchor=north,
        font=\scriptsize,
        nodes={scale=0.8, transform shape},
        legend columns=2,
        /tikz/every even column/.style={column sep=5pt}
    },
    tick label style={font=\small},
    label style={font=\small},
    title style={font=\scriptsize\bfseries},
    xmin=-2, xmax=2
]

\addplot[color=blue!20, thick]
    table[x=x,y=rho_0.1,col sep=comma]
    {tikz/rho-limit-study/parabolas-data.csv};
\addlegendentry{$\rho=10^{-1}$}

\addplot[color=blue!45, thick]
    table[x=x,y=rho_0.01,col sep=comma]
    {tikz/rho-limit-study/parabolas-data.csv};
\addlegendentry{$\rho=10^{-2}$}

\addplot[color=blue!70, thick]
    table[x=x,y=rho_0.001,col sep=comma]
    {tikz/rho-limit-study/parabolas-data.csv};
\addlegendentry{$\rho=10^{-3}$}

\addplot[color=blue!95, thick]
    table[x=x,y=rho_0.0001,col sep=comma]
    {tikz/rho-limit-study/parabolas-data.csv};
\addlegendentry{$\rho=10^{-4}$}

\addplot[color=red, dashed, ultra thick, domain=-2:2, samples=2] {2.616867};
\addlegendentry{Benchmark $\lambda_{\Omega,h}$}

\end{axis}
\end{tikzpicture}
\begin{tikzpicture}
\begin{axis}[
    scale only axis,
    width=6.25cm,
    height=4.5cm,
    xmode=log,
    ymode=log,
    x dir=reverse, 
    title={Log-Log Convergence ($T=10$)},
    xlabel={Discount Factor $\rho$},
    ylabel={Max Error $\|\rho u_\rho - \lambda_{\Omega, h}\|_\infty$},
    grid=both,
    grid style={dashed, gray!30},
    legend pos=north east,
    legend cell align={left},
    legend style={
        at={(0.98,0.98)}, 
        anchor=north east, 
        font=\scriptsize, 
        nodes={scale=0.8, transform shape}
    },
    tick label style={font=\small},
    label style={font=\small},
    title style={font=\scriptsize\bfseries}
]

\addplot[
    color=blue!70!black,
    mark=*,
    thick,
    mark options={fill=blue!70!black}
] table[x=rho,y=error,col sep=comma]
    {tikz/rho-limit-study/convergence-data.csv};
\addlegendentry{Numerical Error}

\addplot[
    color=gray,
    dashed,
    thick
] coordinates {
    (1.0000e-01, 1.2000e-01)
    (1.0000e-04, 1.2000e-04)
};
\addlegendentry{$\mathcal{O}(\rho)$ Reference}

\end{axis}
\end{tikzpicture}    
 \caption{Numerical verification of the discrete vanishing-discount limit $\rho u_\rho \rightarrow \lambda_{\Omega,h}$ in $d=1$ on $\Omega=[-2,2]$, with $\sigma=2$, $\nu=2$, $N_h=200$, $N_T=100$, and $T=10$.
 The stationary discrete benchmark, computed by the Newton method on the same grid, is $\lambda_{\Omega,h}=2.616867$.
 \textit{(Left)} Spatial profiles of $\rho u_\rho$ for $\rho \in \{10^{-1},10^{-2},10^{-3},10^{-4}\}$.
\textit{(Right)} Log-log plot of $\|\rho u_\rho(\cdot,T/2)-\lambda_{\Omega,h}\|_\infty$, consistent with $\mathcal{O}(\rho)$ scaling.}
 \label{fig:rho-limit-study}
\end{figure}

For the sweep reported in Figure~\ref{fig:rho-limit-study}, we use $\Omega=[-2,2]$, $\sigma=2$, $\nu=2$, $N_h=200$, $N_T=100$, and $T=10$.
{We first solve the stationary discrete ergodic system by the Newton method with Dirichlet data from the closed-form profiles, obtaining $\lambda_{\Omega,h}=2.616867$, close to the continuum value $\lambda=2.612086$ from \eqref{eq:true-ergodic-constant}.}
The mass constraint in this stationary solve uses the truncated mass of the closed-form Gaussian profile on $\Omega$, rather than renormalizing it to one, so that the stationary reference and the discounted finite-horizon solve use the same finite-domain window.
We then solve the finite-horizon discounted system \eqref{eq:forward-backward-discounted-PDE-system} for decreasing values of $\rho$, using the Newton method warm-started by continuation in $\rho$, and extract the midpoint product $\rho u_\rho(\cdot,T/2)$.
{The Dirichlet data for this sweep are $\bar m_\rho$ and $\eta_\rho|x|^2+\lambda_{\Omega,h}/\rho$, so that $\rho u_\rho$ equals $\lambda_{\Omega,h}+4\rho\eta_\rho$ on the boundary. Since $\|U\|_\infty=\mathcal{O}(1/\rho)$, the Newton stopping test in this sweep uses the residual scaled by $\max(1,\|U\|_\infty)$ with tolerance $2\times10^{-9}$.}

At fixed $\Omega$, $h$, and $T$, this comparison includes three numerical effects: domain truncation, spatial discretization, and finite-horizon turnpike error.
The first two are shared by the stationary reference and the dynamic discounted solve; the last is assessed empirically through the flattening of $\rho u_\rho(\cdot,T/2)$ and its error relative to $\lambda_{\Omega,h}$.
Thus Figure~\ref{fig:rho-limit-study} should be read as a numerical verification of the discrete vanishing-discount behavior for this grid and horizon, rather than as an independent convergence theorem.

Although we use closed-form solutions for validation here, the same framework can test the vanishing-discount limit when no closed-form solution is available.
In such general cases, the transition from a ``Dirichlet window'' to an arbitrary physical domain is achieved by replacing boundary anchoring with zero-flux Neumann boundaries, treating the grid as a closed physical system. 
We then anchor the stationary solution by imposing the normalization $\int_{\Omega}u\,dx=0$ and apply the corresponding mean subtraction to the dynamic solution.

The results, visualized in Figure~\ref{fig:rho-limit-study}, show the product $\rho u_\rho$ flattening spatially toward the discrete benchmark $\lambda_{\Omega,h}$.
To quantify this behavior, we track $\|\rho u_\rho(\cdot, T/2)-\lambda_{\Omega,h}\|_\infty$.
The resulting log-log plot is consistent with the expected scaling
$\|\rho u_\rho(\cdot,T/2)-\lambda_{\Omega,h}\|_\infty=\mathcal{O}(\rho)$, providing a numerical check of the discrete vanishing-discount relation on the chosen grid.
{The maximum is attained on the boundary, where the prescribed data differ from $\lambda_{\Omega,h}$ by exactly $4\rho\eta_\rho$; the interior values, for instance the maximum over $|x|\le 1$, decrease at the same first-order rate, while the value at the center decreases faster.}

\subsection{A problem with nonlocal congestion effects}
\label{sec:numexp:nonlocalcongestion}

We first introduce a local congestion model with a non-separable Hamiltonian.
In this setting, the cost of motion increases with the local density. 
This introduces a stronger coupling between the HJB and FP equations, since the optimal drift depends explicitly on the mass distribution $m$. 
We solve the problem on the torus $\mathcal{D}=\mathbb{T}^d = [0, 1]^d$. 
The corresponding Hamiltonian is
\begin{equation} \label{eq:local-hamiltonian}
    H(x, m, p) = \frac{|p|^2}{2(1 + 4m)^\beta} - \zeta m.
\end{equation}

The denominator $(1 + 4m)^\beta$ models congestion: as the density $m$ increases, mobility decreases, making rapid movement increasingly costly. 
The term $-\zeta m$ in the Hamiltonian corresponds to a $\zeta m$ contribution in the running cost: agents pay a linear penalty proportional to the local density, discouraging them from remaining in crowded regions. 
This additive crowd-aversion penalty acts alongside the multiplicative mobility reduction $(1+4m)^{\beta}$ in the kinetic prefactor.

We now introduce nonlocal congestion by replacing the pointwise density in~\eqref{eq:local-hamiltonian} with a spatial average. Let
\begin{equation*}
\bar m(x) = (K\ast m)(x) = \int_{\mathcal{D}} K(x-y)\,m(y)\,dy
\end{equation*}
denote the spatially averaged density. The corresponding nonlocal
Hamiltonian, viewed as a functional of the density, is
\begin{equation} \label{eq:nonlocal-hamiltonian}
H_{\mathrm{nl}}[m](x,p)
  = \frac{|p|^2}{2\!\left(1 + 4\bar m(x)\right)^{\!\beta}}
  - \zeta \bar m(x).
\end{equation}
The convolution kernel $K$ is a Gaussian function of the distance on the torus:
\begin{equation} \label{eq:nonlocal-kernel}
  K(z) = \frac{1}{Z}
         \exp\!\left(-\frac{|z|^2}{2\sigma^2}\right),
\end{equation}
where $|z|$ denotes the distance on the torus $\mathcal{D}$ and $Z$ is chosen so that $\int_{\mathcal{D}} K(z)\,dz = 1$.
The standard deviation $\sigma$ governs the radius of the nonlocal interaction.
The discrete convolution is computed by FFT, corresponding to circular convolution with the periodized kernel. The kernel transform is precomputed once and reused at every Picard and Newton step.

The local model is recovered as the sensing radius tends to zero at fixed viscosity. Writing $K_\sigma$ to make the kernel width explicit, $K_\sigma\ast m\to m$ at continuity points of the density as $\sigma\to0$, and the nonlocal Hamiltonian~\eqref{eq:nonlocal-hamiltonian} tends to the local Hamiltonian~\eqref{eq:local-hamiltonian}. Thus local congestion corresponds to the limit of vanishing interaction radius.

For the Picard method, the density iterate is held fixed during each backward HJB sweep, so $\partial\tilde H/\partial M$ does not enter the HJB Newton matrix. The HJB and FP spatial matrices retain their tridiagonal structure in 1D and pentadiagonal structure in 2D, allowing direct sparse solves as described in Section~\ref{sec:numericalmethods:practical}. The nonlocal density enters the residuals and local coefficients through FFT-based convolution, at $O(N_h^2\log N_h)$ cost per timestep in 2D, without assembling dense density-derivative blocks.

The introduction of the nonlocal integral operator fundamentally alters the matrices arising from the Newton method.
In the locally separable and local congestion cases, the derivative $\partial \tilde{H}_i / \partial M_j$ is zero for $i \neq j$, ensuring that the cross-coupling block $A_{\cU,\cM}$ in the global Jacobian remains sparse.
The nonlocal kernel, however, couples $\tilde{H}_i$ to $M_j$ at every node $j$ within the kernel's support.
Because the Gaussian kernel has full support on the torus, the blocks $E_n$ comprising $A_{\cU,\cM}$ are effectively fully dense matrices of size $N_h^d \times N_h^d$, and the off-diagonal blocks of $A_{\cM,\cM}$ inherit the same dense structure through the linearization of $B$ in $M$.

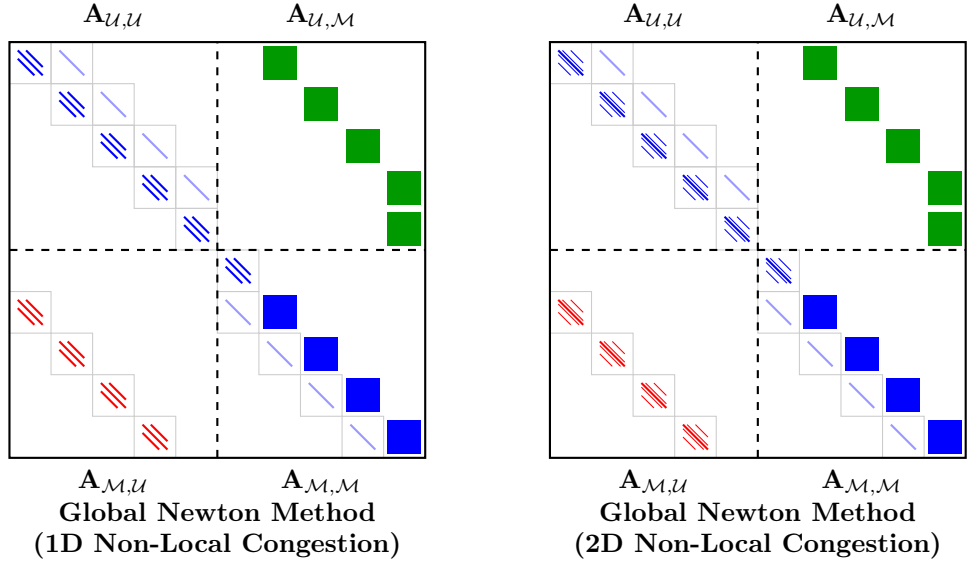
\begin{figure}[H]
    \centering
\newcommand{\tridiagbox}[3]{
    \draw[gray!40, thin] (#1, #2) rectangle (#1+1, #2+1);
    \draw[#3, thick] (#1+0.2, #2+0.8) -- (#1+0.8, #2+0.2);
    \draw[#3, thick] (#1+0.4, #2+0.8) -- (#1+0.8, #2+0.4);
    \draw[#3, thick] (#1+0.2, #2+0.6) -- (#1+0.6, #2+0.2);
}
\newcommand{\pentadiagbox}[3]{
    \draw[gray!40, thin] (#1, #2) rectangle (#1+1, #2+1);
    \draw[#3, thick] (#1+0.2, #2+0.8) -- (#1+0.8, #2+0.2);
    \draw[#3, thin]  (#1+0.3, #2+0.8) -- (#1+0.8, #2+0.3);
    \draw[#3, thin]  (#1+0.2, #2+0.7) -- (#1+0.7, #2+0.2);
    \draw[#3, thin]  (#1+0.5, #2+0.8) -- (#1+0.8, #2+0.5);
    \draw[#3, thin]  (#1+0.2, #2+0.5) -- (#1+0.5, #2+0.2);
}
\newcommand{\densebox}[3]{
    \fill[#3] (#1+0.1, #2+0.1) rectangle (#1+0.9, #2+0.9);
}
\newcommand{\diagbox}[3]{
    \draw[gray!40, thin] (#1, #2) rectangle (#1+1, #2+1);
    \draw[#3, thick] (#1+0.2, #2+0.8) -- (#1+0.8, #2+0.2);
}
\begin{tikzpicture}[scale=0.55]
    \begin{scope}[shift={(0,0)}]
        \foreach \i in {0,1,2,3,4} { \tridiagbox{\i}{9-\i}{blue} }
        \foreach \i in {1,2,3,4}   { \diagbox{\i}{10-\i}{blue!40} }

        \tridiagbox{5}{4}{blue}
        \foreach \i in {1,2,3,4} { \densebox{5+\i}{4-\i}{blue} }
        \foreach \i in {0,1,2,3} { \diagbox{5+\i}{3-\i}{blue!40} }

        \foreach \i in {0,1,2,3} { \tridiagbox{\i}{3-\i}{red} }

        \foreach \i in {1,2,3,4} { \densebox{5+\i}{10-\i}{green!60!black} }
        \densebox{9}{5}{green!60!black}

        \draw[dashed, thick] (5,0) -- (5,10);
        \draw[dashed, thick] (0,5) -- (10,5);
        \draw[black, thick] (0,0) rectangle (10,10);

        \node at (2.5, 10.6) {$\mathbf{A_{\mathcal{U},\mathcal{U}}}$};
        \node at (7.5, 10.6) {$\mathbf{A_{\mathcal{U},\mathcal{M}}}$};
        \node at (2.5, -0.6) {$\mathbf{A_{\mathcal{M},\mathcal{U}}}$};
        \node at (7.5, -0.6) {$\mathbf{A_{\mathcal{M},\mathcal{M}}}$};

        \node at (5, -1.3) {\textbf{Global Newton Method}};
        \node at (5, -2.1) {\textbf{(1D Non-Local Congestion)}};
    \end{scope}

    \begin{scope}[shift={(13,0)}]
        \foreach \i in {0,1,2,3,4} { \pentadiagbox{\i}{9-\i}{blue} }
        \foreach \i in {1,2,3,4}   { \diagbox{\i}{10-\i}{blue!40} }

        \pentadiagbox{5}{4}{blue}
        \foreach \i in {1,2,3,4} { \densebox{5+\i}{4-\i}{blue} }
        \foreach \i in {0,1,2,3} { \diagbox{5+\i}{3-\i}{blue!40} }

        \foreach \i in {0,1,2,3} { \pentadiagbox{\i}{3-\i}{red} }

        \foreach \i in {1,2,3,4} { \densebox{5+\i}{10-\i}{green!60!black} }
        \densebox{9}{5}{green!60!black}

        \draw[dashed, thick] (5,0) -- (5,10);
        \draw[dashed, thick] (0,5) -- (10,5);
        \draw[black, thick] (0,0) rectangle (10,10);

        \node at (2.5, 10.6) {$\mathbf{A_{\mathcal{U},\mathcal{U}}}$};
        \node at (7.5, 10.6) {$\mathbf{A_{\mathcal{U},\mathcal{M}}}$};
        \node at (2.5, -0.6) {$\mathbf{A_{\mathcal{M},\mathcal{U}}}$};
        \node at (7.5, -0.6) {$\mathbf{A_{\mathcal{M},\mathcal{M}}}$};

        \node at (5, -1.3) {\textbf{Global Newton Method}};
        \node at (5, -2.1) {\textbf{(2D Non-Local Congestion)}};
    \end{scope}
\end{tikzpicture}
    \caption{Sparsity structure for a non-separable Hamiltonian with nonlocal congestion in 1D \textit{(Left)} and 2D \textit{(Right)}.
    Each block represents an $N_h \times N_h$ (1D) or $N_h^2 \times N_h^2$ (2D) spatial operator. 
    The cross-coupling blocks $E_n$ in $A_{\mathcal{U},\mathcal{M}}$ consist of dense submatrices (solid green squares) because the nonlocal congestion term couples the value function $U$ to the mass $M$ across the entire spatial domain at each timestep. 
    The off-diagonal blocks in $A_{\mathcal{M},\mathcal{M}}$ (solid blue squares) are likewise dense, since linearizing the transport operator $B$ with respect to $M$ inherits the kernel's full-domain coupling.
    The resulting dense blocks in $A_{\mathcal{U},\mathcal{M}}$ and $A_{\mathcal{M},\mathcal{M}}$ sharply increase the memory and matrix-vector product cost of the assembled Jacobian relative to the locally separable case.}
\label{fig:sparsity_nonseparable_nonlocal_H_1d2d}
\end{figure}

In $d=1$ this transition remains computationally manageable.
For $N_h = 200$, each dense block $E_n$ contains $40{,}000$ entries; across $N_T = 100$ timesteps this is a memory overhead of approximately $32\,\text{MB}$, which keeps the assembled space-time Jacobian well within reach of the assembled direct sparse backend in one dimension.

In $d=2$, scaling the Newton method to the nonlocal regime under our common linear-algebra backend becomes intractable.
For $N_x = N_y = 50$, a single dense $E_n$ already contains $(2{,}500)^2 = 6.25 \times 10^6$ entries; at the target resolution $N_x = N_y = 200$ this rises to $(40{,}000)^2 = 1.6 \times 10^9$ entries per block, so assembling the space-time Jacobian over $N_T = 20$ timesteps would require in excess of $250\,\text{GB}$ of RAM for the $A_{\cU,\cM}$ quadrant alone.
Arithmetic on the assembled dense blocks scales as $O(N_T \cdot N_h^4)$ per Jacobian application, compounding the memory barrier.

This limitation is specific to the assembled Newton backend used here.
A matrix-free implementation could exploit the convolution structure of the nonlocal terms to apply the dense blocks without storing them explicitly.

Within the assembled-matrix framework considered here, this difference favors the Picard method for the two-dimensional nonlocal problem.

Both models use the following initial densities and terminal costs.
For $d=1$, the system is initialized with a ``square-wave'' mass distribution concentrated in the center of the domain:
\begin{equation}m_0(x) = C \mathds{1}_{[0.375, 0.625]}(x)
\end{equation}
where $C$ is chosen so that $\int m_0 \ dx = 1$.
The terminal cost is a double quadratic well that encourages agents to split and move toward two targets located at $x=0.3$ and $x=0.7$:
\begin{equation}g(x) = 15 \min\left( (x - 0.3)^2, (x - 0.7)^2 \right).\end{equation}
This setup creates a conflict: agents want to move rapidly from the center to the targets to minimize the terminal cost, but the high initial density in the center combined with the congestion term makes this initial movement expensive.
For $d=2$, we naturally extend the setup to the square torus. The initial mass is defined as a square patch of uniform density in the center of the domain:
\begin{equation}
m_0(x, y) = C \mathds{1}_{[0.375, 0.625]}(x) \mathds{1}_{[0.375, 0.625]}(y).
\end{equation}
We extend the terminal cost by placing two symmetric wells along the $x$-axis, thereby reproducing the splitting dynamics of the 1D problem:
\begin{equation}
g(x, y) = 15 \min\left( (x - 0.3)^2 + (y-0.5)^2, (x - 0.7)^2 + (y-0.5)^2 \right).
\end{equation}

\begin{figure}[H]
 \centering
\includegraphics[width=1\textwidth]{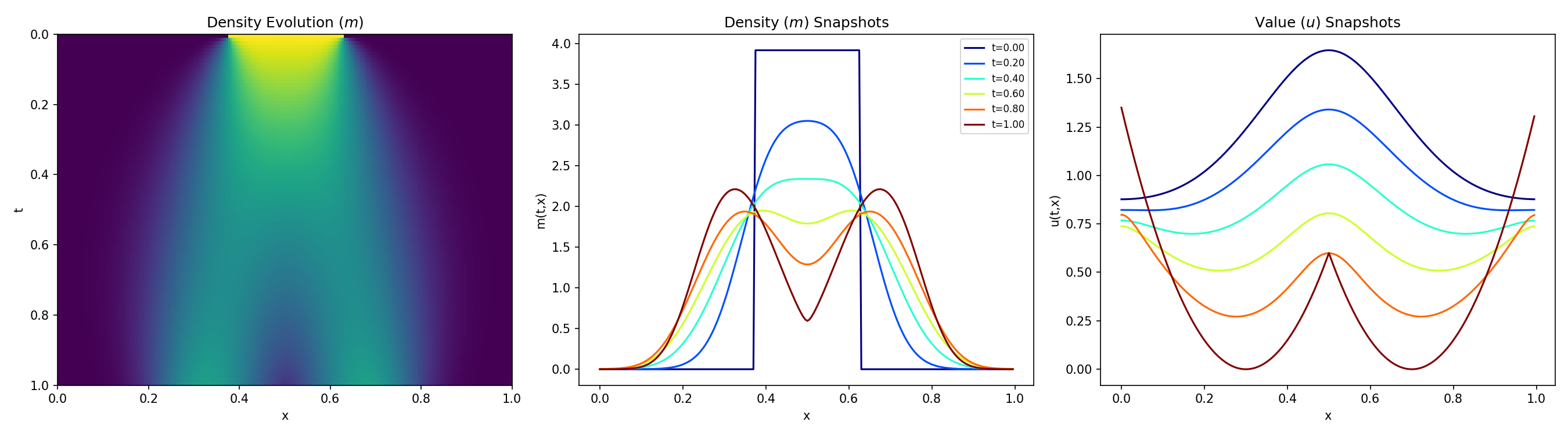}
\caption{The nonlocal congestion problem in $d = 1$.
\textit{(Left)} Density heatmap showing a more gradual dispersal from the initial square-wave than in the local case.
\textit{(Middle)} Density snapshots $m(t,x)$: the nonlocal profile maintains a higher central peak before bifurcation, and the terminal peaks at $T=1$ reach higher amplitudes than in the local model.
\textit{(Right)} Value function snapshots $u(t,x)$, smoother throughout the time horizon than in the local case as a result of the nonlocal operator.} 
 \label{fig:nonlocal-congestion-slice-1d}
\end{figure}

\begin{figure}[H]
 \centering
 \includegraphics[width=1\textwidth]{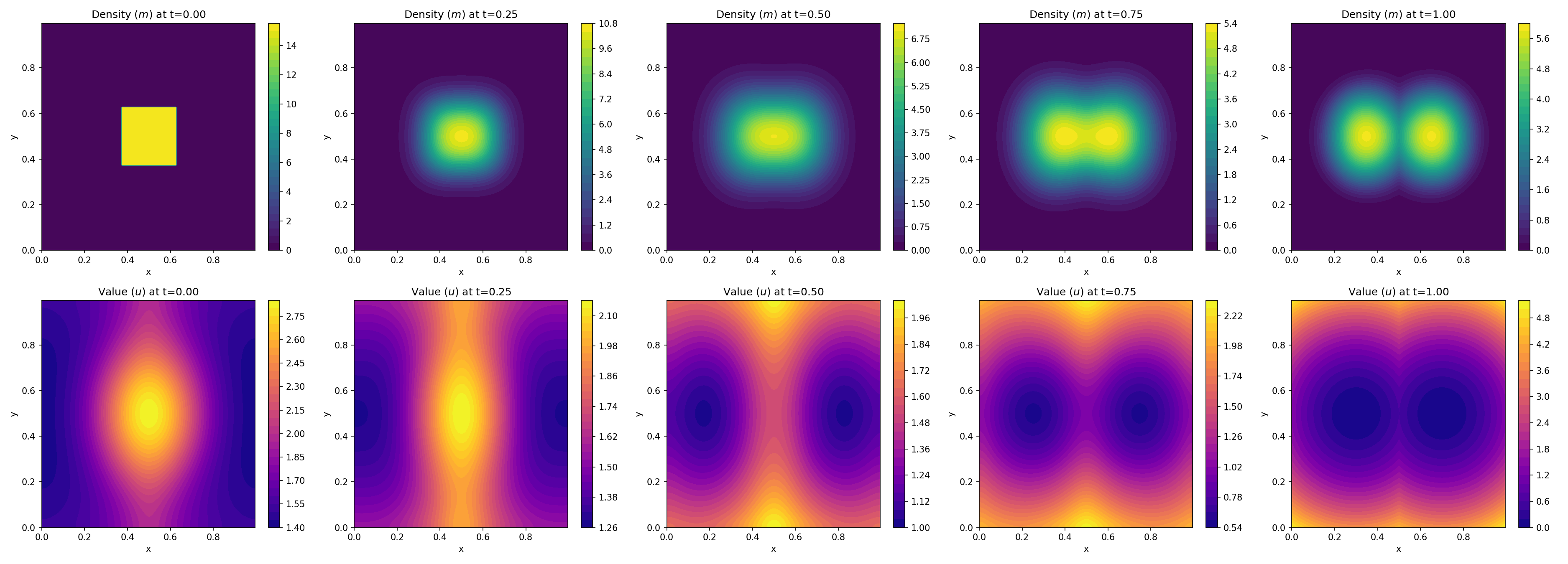}
\caption{The nonlocal congestion problem in $d = 2$. 
\textit{(Top)} Density contour plots (left to right, $t = 0, 0.25, 0.5, 0.75, 1.0$) showing the initial square patch maintaining a concentrated single-hump profile, followed by slow symmetric radial dispersal under the attenuated congestion, and finally two sharp terminal peaks with amplitudes substantially higher than in the local 2D solution.
\textit{(Bottom)} Value function evolution from a broad, smooth central bowl to the double-well terminal structure, with smoother contours throughout than in the local case.}
 \label{fig:nonlocal-congestion-slice-2d}
\end{figure}

We set $\beta = 1.5$, $\zeta = 1$, $\nu = 0.01$, $T=1$, and the sensing radius
$\sigma = 0.2$.  The nonlocal Picard runs use $N_h=200$ in each coordinate direction, with $N_T=100$ in $d=1$ and $N_T=20$ in $d=2$; both use $\theta=0.2$.
The nonlocal coupling introduces significantly stiffer Picard dynamics than the local case.
The nonlocal problem uses $\theta = 0.2$ to prevent the Picard iterates from diverging.
With $\sigma = 0.2$ on a unit domain, the kernel's effective support (roughly $\pm 3\sigma$) spans most of the torus, so a perturbation to the mass anywhere in the domain immediately perturbs the Hamiltonian almost everywhere, amplifying the sensitivity between successive iterates.
Figures~\ref{fig:nonlocal-congestion-slice-1d} and~\ref{fig:nonlocal-congestion-slice-2d} show the evolution of the density $m$ and value function $u$ in $d = 1$ and $d = 2$, respectively.

The key to understanding the nonlocal case is that agents react not to the pointwise density $m$ but to the spatially averaged density $\bar m(x) = (K \ast m)(x)$.
Because the Gaussian kernel has effective support over several standard deviations, the convolution averages the density over a region much larger than the initial square-wave support.
Across the initial blob, $\bar m$ is therefore lower than $m$, because the kernel integrates over the empty space surrounding the blob.
Consequently, the nonlocal distribution remains more concentrated than the local distribution throughout the time horizon.
Weaker dispersive pressure means that agents have less incentive to flee crowded regions, so they linger longer at the center before bifurcating.
At the terminal wells, the effective congestion cost felt by each agent is therefore reduced, lowering the penalty for tight concentration there.
This behavior is consistent with Aurell and Djehiche~\cite{Aurell2018nonlocalcrowd}, who observe that widening the interaction radius causes the equilibrium to transition from a strongly dispersed, local-like regime to a regime in which agents can concentrate more tightly because each perceives only the averaged, attenuated congestion of its neighbors.

The two-dimensional results confirm the same picture.
The initial square patch disperses more slowly than in the local 2D case; the bifurcation does not become pronounced until $t\approx 0.75$.
At the terminal time, the density peaks at the two wells are higher than those in the local two-dimensional solution.

\subsection{A ``traffic-light'' problem: an example of a temporal ``shock''}
\label{sec:numexp:traffic-light}

This benchmark features a time-dependent obstacle that produces a temporal shock: agents traverse a one-dimensional domain from left to right but must navigate a ``traffic light'' at its center.
This light starts green, turns red for a specified interval in the middle of the time horizon, and subsequently turns green again. 
We use the Hamiltonian
$$H(x, m, p) = \frac{1}{2}|p|^2 - V(x, t) -\kappa m^{\alpha}.$$

The term $m^{\alpha}$ represents a congestion cost, penalizing agents for moving through areas of high density. This reflects agents' tendency to avoid crowds and the greater effort required to navigate heavy traffic.

We model the traffic light with a time-dependent potential
\begin{equation}
    V(x,t) = P\,I(t)\,G_w(x;\mu), 
    \qquad 
    G_w(x;\mu) = \exp\!\left(-\frac{(x-\mu)^2}{2(w/2)^2}\right).
\end{equation}

The light is represented by a steep Gaussian centered at $\mu$, the midpoint of the spatial domain, with width $w$. 
The indicator $I(t)$ is active on the red-light interval $[\delta T, (1-\delta)T]$, where $\delta = 0.3$, and $P$ is a large penalty for running the red light.
The spatial Gaussian is a modeling choice for the obstacle's shape and plays no role in the numerical difficulty. The shock is temporal: $I(t)$ is a step function, and $V(x, t)$ changes abruptly between consecutive timesteps at $t = 0.3$ and $t = 0.7$.

For this experiment, we take $\mathcal{D} = [0, 1]$ and $T = 1$. 
Homogeneous Neumann boundary conditions model a closed corridor: agents accumulate at $x = 1$ rather than exiting the domain.
The initial mass is a Gaussian centered at $x = 0.1$,
\begin{equation}
m_0(x) \propto \exp\!\left(-\frac{(x - 0.1)^2}{0.01}\right),
\end{equation}
normalized so that $\int_{\mathcal{D}} m_0(x)\,dx = 1$.
The terminal cost, $g(x) = 5(x-1)^2$, strongly compels agents to reach the right endpoint ($x = 1$) of the domain.
The traffic light obstacle is a Gaussian of width $w = 0.05$ centered at $\mu = 0.5$, with penalty $P = 150$.

The light remains active from $t = 0.3$ to $t = 0.7$, yielding the indicator function

\[ I(t) = \begin{cases} 
      1 & 0.3 \leq t \leq 0.7 \\
      0 & \text{elsewhere}. 
   \end{cases}
\]

We set the congestion exponent to $\alpha= 2$, the congestion coefficient to $\kappa = 1$, and the viscosity coefficient to $\nu = 0.02$. 
The discretization is fixed at $N_h = 200$ and $N_T = 100$.

\begin{figure}[H]
    \centering
    \includegraphics[width=0.98\textwidth]{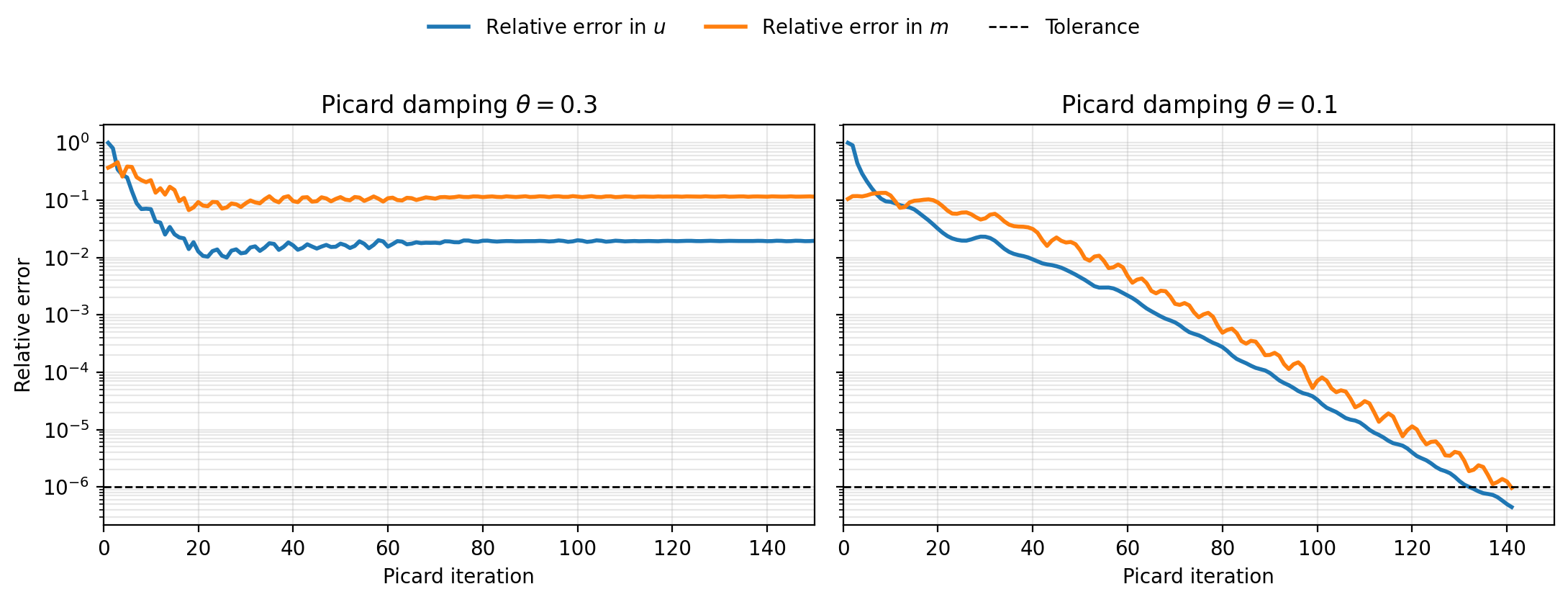}
    \caption{Picard iteration history for the traffic-light problem, demonstrating convergence failure for $\theta = 0.3$ \textit{(Left)} and convergence in 141 iterations for $\theta = 0.1$ \textit{(Right)}.}
    \label{fig:trafficlight-picard-history}
\end{figure}

\begin{figure}[p]
 \centering
 \includegraphics[width=0.475\textwidth]{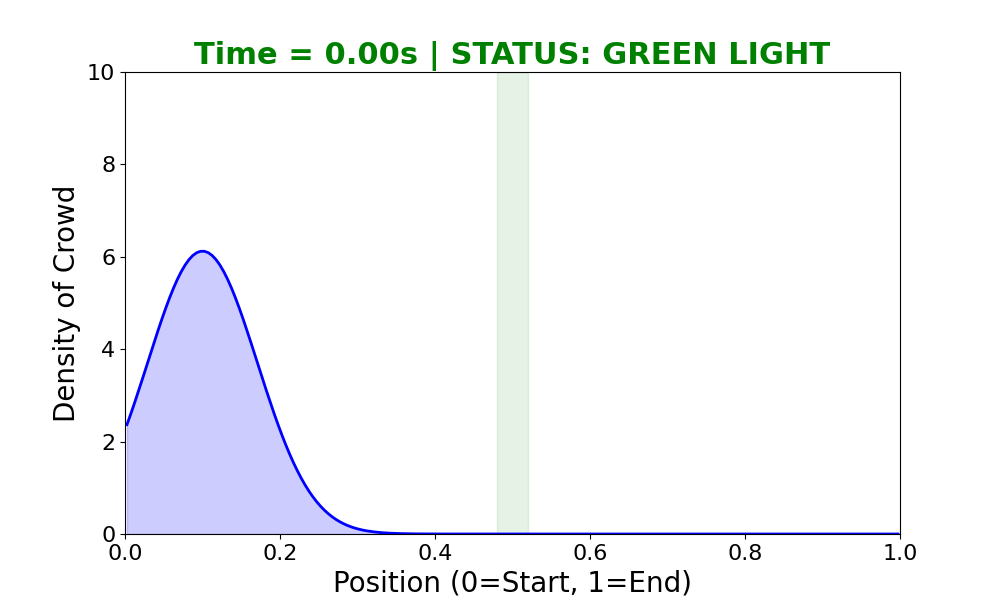}
 \includegraphics[width=0.475\textwidth]{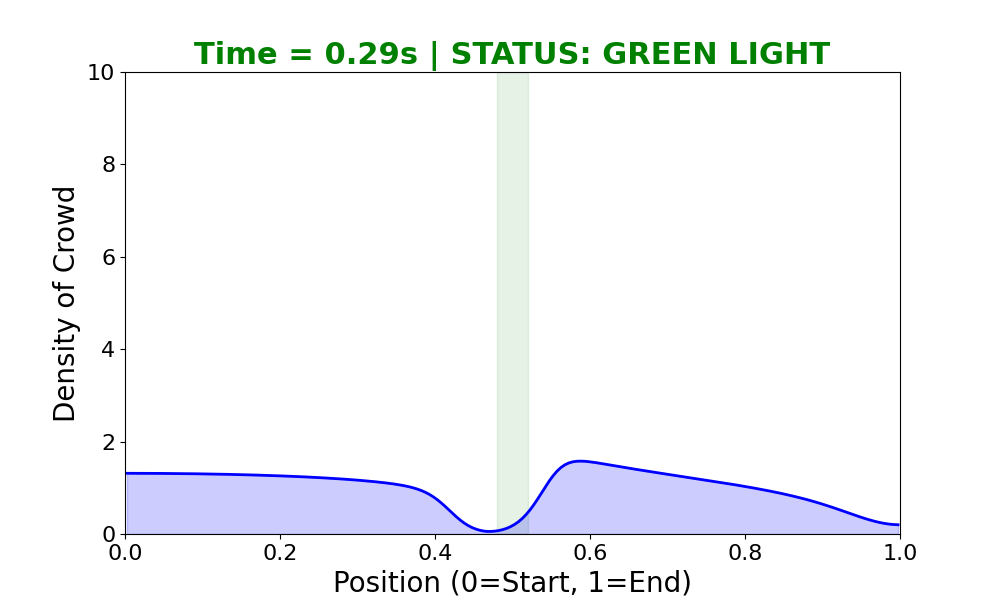}
 \includegraphics[width=0.475\textwidth]{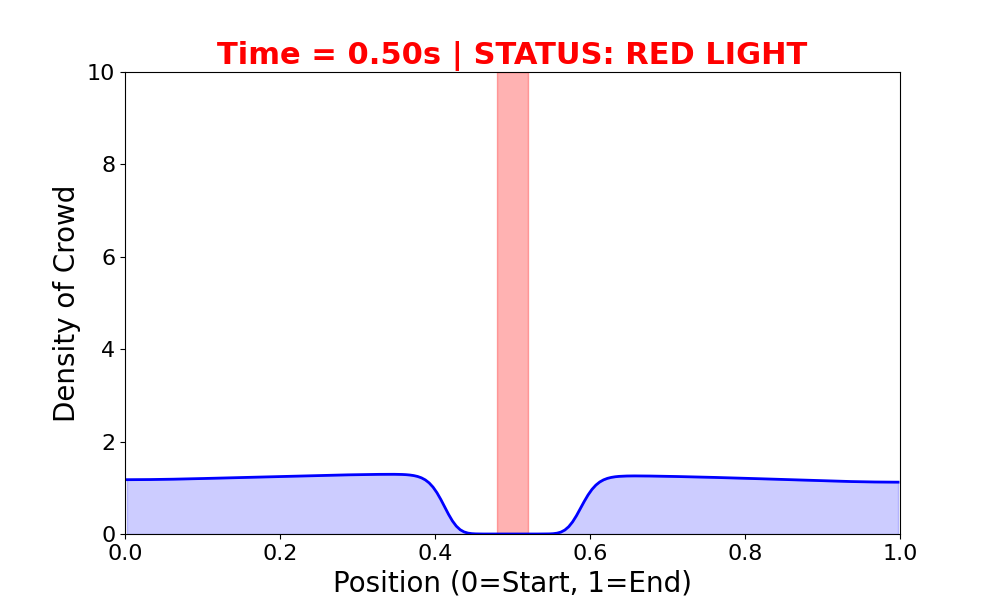}
 \includegraphics[width=0.475\textwidth]{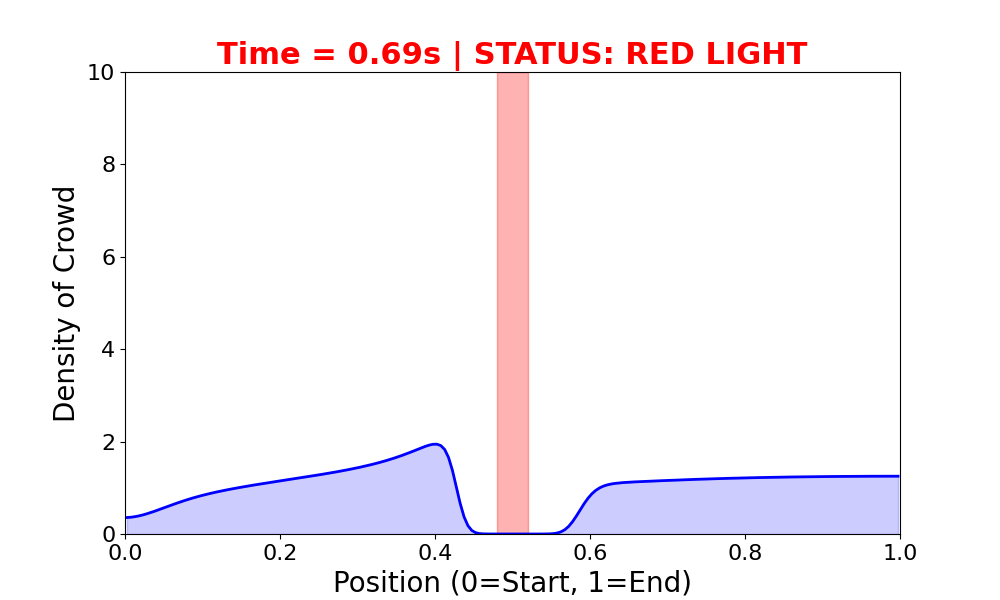}
 \includegraphics[width=0.475\textwidth]{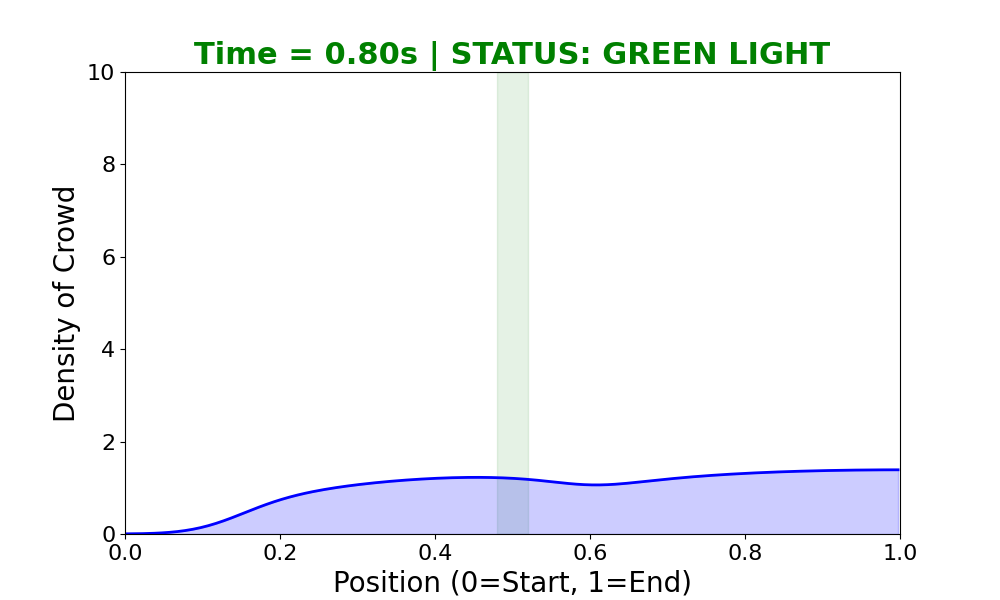}
 \includegraphics[width=0.475\textwidth]{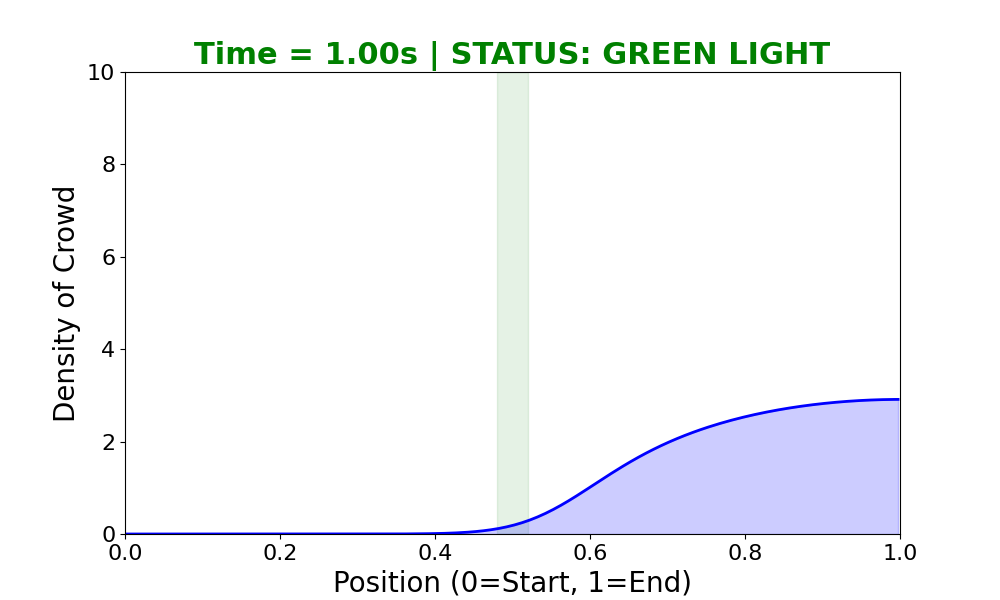}
 \caption{Snapshots of the density in the traffic-light problem. \textit{(Top left)} Initial distribution at $t = 0$. \textit{(Top right)} Distribution at $t = 0.29$, showing agents anticipating the red-light interval beginning at $t = 0.3$. \textit{(Middle left)} Distribution at the midpoint of the red-light interval. \textit{(Middle right)} Distribution at $t = 0.69$, showing agents anticipating the end of the red-light interval at $t = 0.7$. \textit{(Bottom left)} Distribution at $t = 0.8$, after the light-induced bifurcation has begun to dissipate. \textit{(Bottom right)} Terminal distribution at $t = 1$, with agents concentrated near the right boundary. Homogeneous Neumann boundary conditions prevent agents from leaving the domain.} 
 \label{fig:trafficlight-picard-snapshots}
\end{figure}

We first solve the system using the Picard method. 
The presence of the temporal shock necessitates significant damping.
As shown in Figure~\ref{fig:trafficlight-picard-history}, the Picard method fails with $\theta = 0.3$ (left) but converges with $\theta = 0.1$ in 141 iterations (right).
As is typical when applying aggressive damping, the Picard loop requires considerably more iterations to achieve convergence.
Figure~\ref{fig:trafficlight-picard-snapshots} shows snapshots of the mass distribution at several times for $\theta = 0.1$.

However, the convergence of the Picard method in this regime remains brittle. 
This problem can be made more difficult by decreasing the viscosity $\nu$ or increasing the congestion exponent $\alpha$.
Newton continuation also has a lower recorded wall-clock cost in this benchmark.
With direct sparse solves on the same grid, $N_h=200$ and $N_T=100$, Newton required $7.55$ seconds including all continuation steps, compared with $25.74$ seconds for the archived Picard run with $\theta=0.1$.
Furthermore, Newton reached $\nu=0.001$ through the aggressive viscosity schedule $0.2\to0.01\to0.001$, whereas the timed Picard run was performed at $\nu=0.02$.
Thus Newton reached a far smaller viscosity in less recorded time.
Newton continuation is also robust to the temporal shock.
Once the continuation process is initiated at a sufficiently high viscosity, the sequence proceeds reliably without stalling. 
The higher-viscosity solution already captures the temporal shock and therefore provides a well-conditioned initial guess for subsequent lower-viscosity solves. This continuation step insulates the Newton iterations from the abrupt traffic-light dynamics.

\FloatBarrier
\subsection{A problem with a double-well potential}
\label{sec:numexp:double-well}

The next benchmark is a finite-horizon problem in which agents start near the center of the domain and move through a tilted double-well potential toward one of its two minima.
The Hamiltonian is
\begin{equation}
H(x, m, p) = \frac{1}{2}|p|^2 - V(x) - \kappa m^{\alpha} \quad (\alpha \ge 1).
\end{equation}

\begin{figure}[H]
 \centering
 \includegraphics[width=1\textwidth]{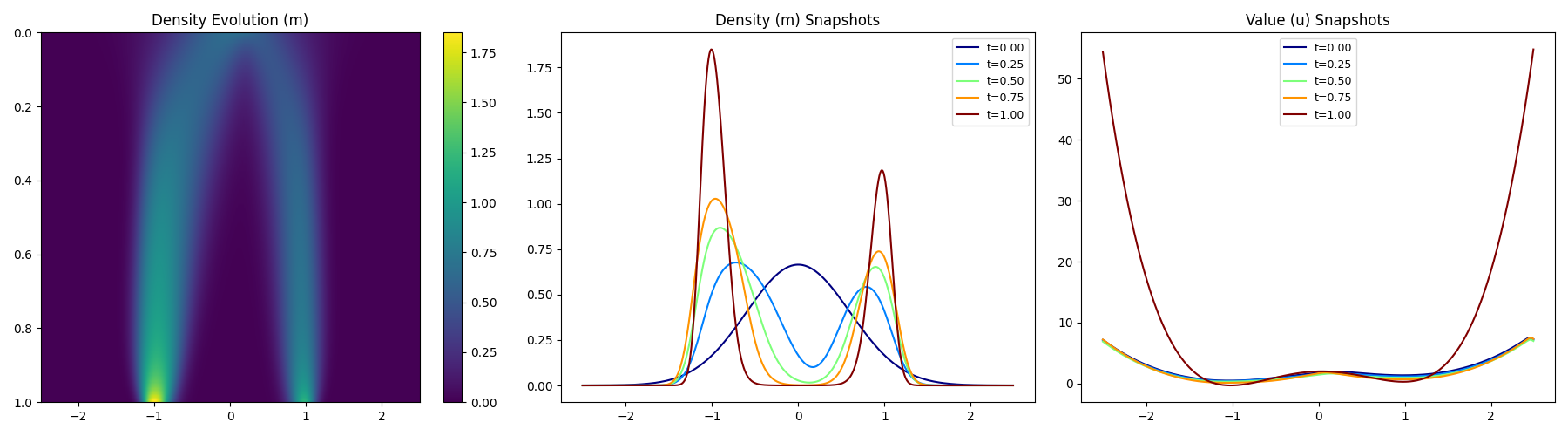}
 \caption{The asymmetric double-well problem in $d = 1$. \textit{(Left)} Density heatmap showing the bifurcation of mass toward two unequal targets. \textit{(Middle)} Evolution of $m(t,x)$ illustrating the preferential concentration in the deeper left well. 
 \textit{(Right)} Value function snapshots showing the transition from a smooth initial state to the terminal double-well potential.} 
 \label{fig:doublewell-slice-1d}
\end{figure}

This corresponds to a running cost $f(x, m, a) = \tfrac{1}{2}|a|^2 + V(x) + \kappa m^{\alpha}$, where $a$ denotes the control; the agents pay a control cost, a potential penalty $V$, and a congestion penalty growing as $m^{\alpha}$.
We take the domain $\Omega=(-2.5,2.5)^d$, centered at the origin, with periodic boundary conditions for both equations, and set $T=1$. 
For $d = 1$, we use the tilted quartic potential
\begin{equation}
V(x) = D(x^2 - 1)^2 + \chi x,
\end{equation}
where $D$ controls the depth of the wells and $\chi$ serves as a linear tilt to break the symmetry. 
Specifically, $\chi > 0$ makes the left well at $x = -1$ lower than the right well at $x = 1$, with $V(-1) - V(1) = -2\chi < 0$, so the left well carries lower terminal cost.
The parameter $\kappa > 0$ serves as a congestion coefficient; because the cost of occupying a state grows as $\kappa m^{\alpha}$, this term forces agents to spread out. 
Without this congestion penalty, all agents would attempt to concentrate exactly at the minima of the wells.
This system is initialized with a Gaussian mass distribution centered at the origin with standard deviation $0.6$, near the local maximum of the tilted potential $V(x)$.
We set the terminal cost to $g(x) = V(x)$ to penalize agents who do not reach the wells by the final time. 
We set $D=2$, $\chi=0.3$, $\alpha=2$, $\kappa=1.0$, and $\nu=0.1$, with space-time resolution $N_h=500$ and $N_T=100$.
Figure~\ref{fig:doublewell-slice-1d} shows the Picard solution in $d = 1$, including a space-time plot of the density evolution and snapshots of the density and value function.

\subsubsection{Picard contractivity diagnostic}

\begin{figure}[H]
    \centering
\begin{tikzpicture}
\begin{axis}[
    scale only axis,
    width=6.75cm,
    height=4.5cm,
    xlabel={$\nu$},
    ylabel={$\rho(P)$},
    x dir=reverse,
    xmin=0.075, xmax=0.315,
    ymin=0.0, ymax=2.25,
xtick={0.10, 0.15, 0.20, 0.30},
xticklabels={$0.10$, $0.15$, $0.20$, $0.30$},
    ytick={0, 0.5, 1.0, 1.5, 2.0},
    grid=major,
    grid style={dashed, gray!60},
    tick label style={font=\small},
    label style={font=\small},
]
\addplot[color=black!40, thin, line width=0.8pt] coordinates {
    (0.30, 0.2936) (0.20, 0.4527) (0.17, 0.6112) (0.15, 0.7646)
    (0.12, 1.1139) (0.10, 1.5561) (0.09, 1.8320) (0.08, 2.2007)
};
\addplot[color=red, dashed, line width=1.5pt] coordinates {
    (0.075, 1.0) (0.315, 1.0)
};
\addplot[color=blue!70!black, only marks, mark=*, mark size=3pt,
    mark options={fill=blue!70!black, draw=black, line width=0.4pt}
] coordinates {
    (0.30, 0.2936) (0.20, 0.4527) (0.17, 0.6112) (0.15, 0.7646)
    (0.12, 1.1139) (0.10, 1.5561) (0.09, 1.8320) (0.08, 2.2007)
};
\end{axis}
\end{tikzpicture}
\begin{tikzpicture}
\begin{axis}[
    scale only axis,
    width=6.75cm,
    height=4.5cm,
    xlabel={$\nu$},
    ylabel={$\rho(P_{\theta})$},
    x dir=reverse,
    xmin=0.075, xmax=0.315,
    ymin=0.0, ymax=2.25,
xtick={0.10, 0.15, 0.20, 0.30},
xticklabels={$0.10$, $0.15$, $0.20$, $0.30$},
    ytick={0, 0.5, 1.0, 1.5, 2.0},
    grid=major,
    grid style={dashed, gray!60},
    tick label style={font=\small},
    label style={font=\small},
]
\addplot[color=black!40, thin, line width=0.8pt] coordinates {
    (0.30, 0.2238) (0.20, 0.2154) (0.17, 0.2882) (0.15, 0.4237)
    (0.12, 0.6966) (0.10, 1.0450) (0.09, 1.2657) (0.08, 1.5607)
};
\addplot[color=red, dashed, line width=1.5pt] coordinates {
    (0.075, 1.0) (0.315, 1.0)
};
\addplot[color=blue!70!black, only marks, mark=*, mark size=3pt,
    mark options={fill=blue!70!black, draw=black, line width=0.4pt}
] coordinates {
    (0.30, 0.2238) (0.20, 0.2154) (0.17, 0.2882) (0.15, 0.4237)
    (0.12, 0.6966) (0.10, 1.0450) (0.09, 1.2657) (0.08, 1.5607)
};
\end{axis}
\end{tikzpicture}
 \caption{Spectral radius of the linearized Picard operator $\rho(P)$ \textit{(Left)} and damped operator $\rho(P_{\theta})$ with $\theta = 0.8$ \textit{(Right)} as a function of viscosity $\nu$, computed for the $d = 1$ double-well problem on a grid with $N_h = 100$ and $N_T = 30$.} 
 \label{fig:spectral-radius-doublewell1d}
\end{figure}
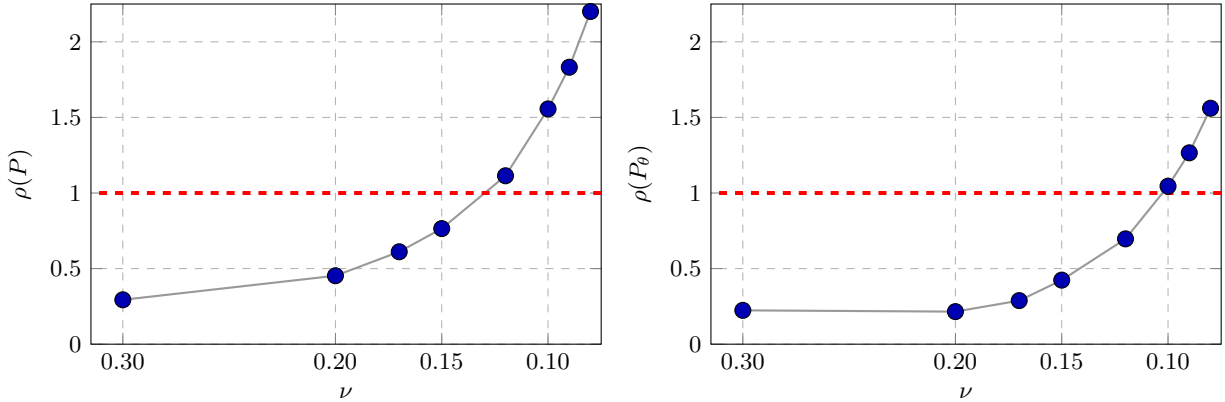

We next use a local spectral-radius diagnostic to explain why the Picard iteration becomes difficult at low viscosity. The upwind finite-difference discretization is designed to remain meaningful in the deterministic limit $\nu \to 0$ \cite{achdou2010numerical}; the issue here is algorithmic rather than a change in the target discrete problem. We linearize the Picard map at a computed fixed point $(U^*,M^*)$ and study the error iteration $\delta U^{(k+1)} = P\delta U^{(k)}$. When $\rho(P)$ crosses 1, the local contraction condition for the Picard iteration is lost.
To derive $P$, we linearize each substep of the Picard iteration around the fixed 
point $(U^*, M^*)$. Let $\delta U^{(k)} = U^{(k)} - U^*$ and $\delta M^{(k)} = 
M^{(k)} - M^*$ denote the errors at iteration $k$.
Given $U^{(k)}$, the FP equation is solved for 
$M^{(k+1)}$. Linearizing $\varphi_M(U^{(k)}, M^{(k+1)}) = 0$ around $(U^*, M^*)$ gives
\begin{equation}
    A_{\cM,\cM}\, \delta M^{(k+1)} = -A_{\cM,\cU}\, \delta U^{(k)}.
\end{equation}
Given $M^{(k+1)}$, the HJB equation is solved 
for $U^{(k+1)}$. Linearizing $\varphi_U(U^{(k+1)}, M^{(k+1)}) = 0$ gives
\begin{equation}
    A_{\cU,\cU}\, \delta U^{(k+1)} = -A_{\cU,\cM}\, \delta M^{(k+1)}.
\end{equation}
Substituting the first relation into the second yields $\delta U^{(k+1)} = 
P\, \delta U^{(k)}$, where
\begin{equation}
    P = A_{\cU,\cU}^{-1} A_{\cU,\cM} A_{\cM,\cM}^{-1} A_{\cM,\cU},
\end{equation}
where the four matrices are the blocks of the Newton Jacobian evaluated at $(U^*,M^*)$.
With the damping parameter $\theta$, the damped iteration satisfies
\begin{equation}
P_{\theta} = (1-\theta)I + 
\theta P.
\end{equation}
For analytical tractability, we restrict the analysis to the variant in which damping is applied only to the $U$ update; the coupled two-variable damping in Algorithm~\ref{alg:picard} admits a similar analysis with the same qualitative conclusion.
By standard linear convergence theory, the damped Picard iteration contracts locally if and only if $\rho(P_{\theta}) < 1$ and the convergence rate is governed by $\rho(P_{\theta})$.
We estimate $\rho(P)$ and $\rho(P_{\theta})$ numerically by power iteration on a finite-difference linearization of the Picard map, on a coarser grid ($N_h = 100$, $N_T = 30$) for tractability; the estimates are indicative of the trend in $\nu$ rather than sharp thresholds.
Figure~\ref{fig:spectral-radius-doublewell1d} reports $\rho(P)$ and $\rho(P_{\theta})$ as a function of $\nu$ for the 1D double-well problem. 
The linearized Picard operator loses its contraction property as $\nu \rightarrow 0$, providing a local explanation for the expected convergence degradation. 
We compute the spectral radius in 1D for tractability; the same mechanism is expected to drive the 2D convergence failure reported below, though we do not compute $\rho(P)$ on the 2D space-time Jacobian directly.
Proposition 1 of \cite{Li2021multiscale} establishes the same local contraction criterion in their alternating-sweeping framework --- with their relaxation factor playing the role of our damping parameter $\theta$ --- though they do not study how the spectral radius depends on $\nu$.

\subsubsection{Two-dimensional cost and hybrid continuation}

\begin{figure}[H]
 \centering
 \includegraphics[width=1\textwidth]{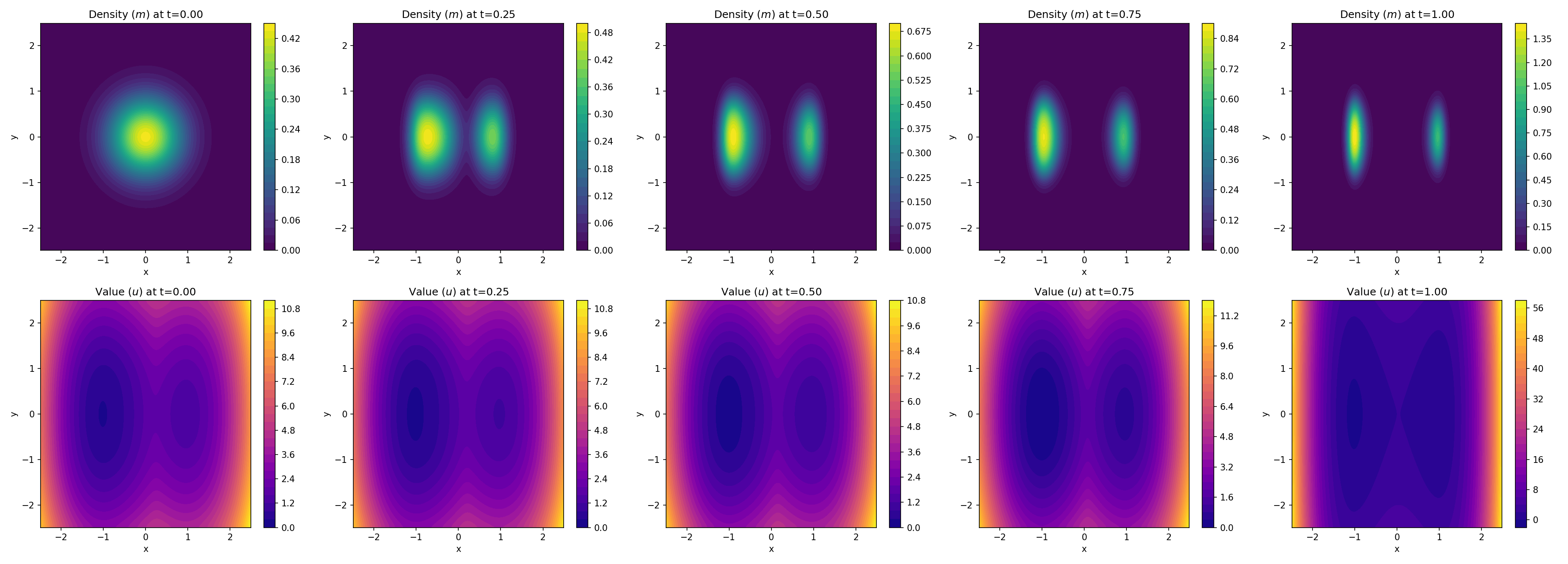}
 \caption{The asymmetric double-well problem in $d = 2$. \textit{(Top)} Density contour plots showing the bifurcation of mass toward the two wells, with preferential concentration in the deeper left well.
 \textit{(Bottom)} Value function evolution toward the terminal double-well potential.} 
 \label{fig:doublewell-slice-2d}
\end{figure}

For the two-dimensional case, we extend the potential to
\begin{equation}
V(x, y) = \left[ D(x^2 - 1)^2 + \chi x \right] + \frac{1}{2} y^2,
\end{equation}
keeping all physical parameters identical.
This potential has two distinct local minima located near $(x, y) = (\pm 1, 0)$.
The addition of the $\frac{1}{2}y^2$ term heavily penalizes deviations from the $x$-axis, effectively focusing the population dynamics along the line connecting the two wells (Figure~\ref{fig:doublewell-slice-2d}). 
The two-dimensional double-well experiment demonstrates how Picard-based acceleration reduces the total cost of reaching a low-viscosity solution: inexpensive Picard solves replace the initial stages of Newton continuation, while Newton completes the descent beyond Picard's observed convergence limit.
For $d=2$, we again use the Picard method and plot the density and value function at several times. 
Relative to the one-dimensional case, only the spatial and temporal resolutions are reduced, to $N_h = 200$ and $N_T = 20$.

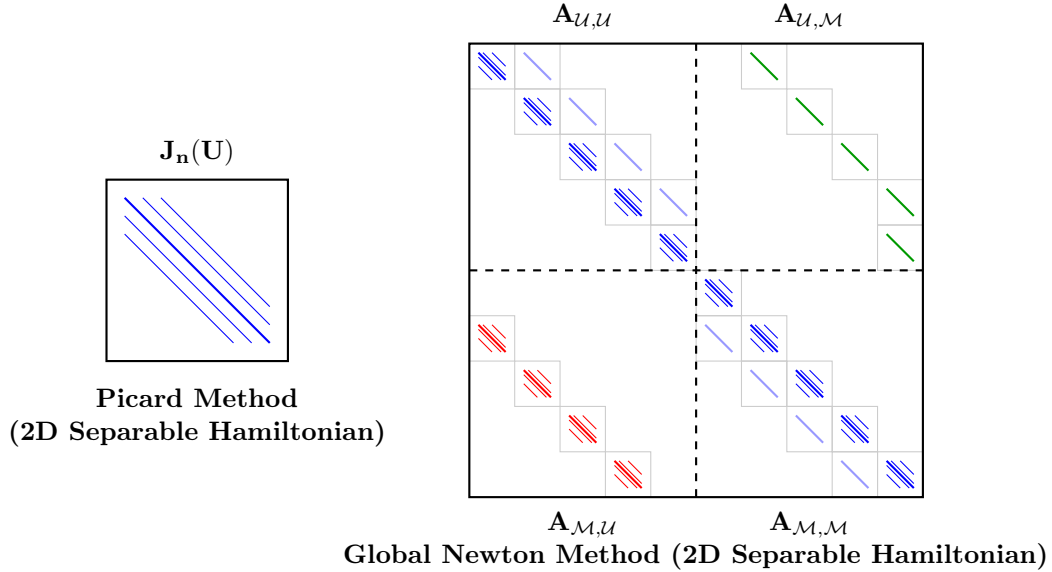
\begin{figure}[H]
\centering
\newcommand{\pentadiagbox}[3]{
    \draw[gray!40, thin] (#1, #2) rectangle (#1+1, #2+1);
    \draw[#3, thick] (#1+0.2, #2+0.8) -- (#1+0.8, #2+0.2);
    \draw[#3, thin] (#1+0.3, #2+0.8) -- (#1+0.8, #2+0.3);
    \draw[#3, thin] (#1+0.2, #2+0.7) -- (#1+0.7, #2+0.2);
    \draw[#3, thin] (#1+0.5, #2+0.8) -- (#1+0.8, #2+0.5);
    \draw[#3, thin] (#1+0.2, #2+0.5) -- (#1+0.5, #2+0.2);
}

\newcommand{\diagbox}[3]{
    \draw[gray!40, thin] (#1, #2) rectangle (#1+1, #2+1);
    \draw[#3, thick] (#1+0.2, #2+0.8) -- (#1+0.8, #2+0.2);
}

\begin{tikzpicture}[scale=0.6]
    \begin{scope}[shift={(0,3)}] 
        \draw[blue, thick] (0.4,3.6) -- (3.6,0.4); 
        \draw[blue, thin] (0.8,3.6) -- (3.6,0.8);  
        \draw[blue, thin] (0.4,3.2) -- (3.2,0.4);
        \draw[blue, thin] (1.2,3.6) -- (3.6,1.2);
        \draw[blue, thin] (0.4,2.8) -- (2.8,0.4);
        
        \draw[black, thick] (0,0) rectangle (4,4);
        \node at (2, 4.6) {$\mathbf{J_n(U)}$};
        
        \node at (2, -0.8) {\textbf{Picard Method}};
        \node at (2, -1.6) {\textbf{(2D Separable Hamiltonian)}};
    \end{scope}

    \begin{scope}[shift={(8,0)}] 
        \foreach \i in {0,1,2,3,4} { \pentadiagbox{\i}{9-\i}{blue} }
        \foreach \i in {1,2,3,4} { \diagbox{\i}{10-\i}{blue!40} }

        \foreach \i in {0,1,2,3,4} { \pentadiagbox{5+\i}{4-\i}{blue} }
        \foreach \i in {0,1,2,3} { \diagbox{5+\i}{3-\i}{blue!40} }

        \foreach \i in {0,1,2,3} { \pentadiagbox{\i}{3-\i}{red} }

        \foreach \i in {1,2,3,4} { \diagbox{5+\i}{10-\i}{green!60!black} }
        \diagbox{9}{5}{green!60!black}

        \draw[dashed, thick] (5,0) -- (5,10);
        \draw[dashed, thick] (0,5) -- (10,5);
        \draw[black, thick] (0,0) rectangle (10,10);
        
        \node at (2.5, 10.6) {$\mathbf{A_{\mathcal{U},\mathcal{U}}}$};
        \node at (7.5, 10.6) {$\mathbf{A_{\mathcal{U},\mathcal{M}}}$};
        \node at (2.5, -0.6) {$\mathbf{A_{\mathcal{M},\mathcal{U}}}$};
        \node at (7.5, -0.6) {$\mathbf{A_{\mathcal{M},\mathcal{M}}}$};
        
        \node at (5, -1.3) {\textbf{Global Newton Method (2D Separable Hamiltonian)}};
    \end{scope}
\end{tikzpicture}
\caption{Sparsity structure for the 2D separable Hamiltonian.
\textit{(Left)} The Picard-method Jacobian $J_n(U)$, which is pentadiagonal of size $N_h^2 \times N_h^2$. \textit{(Right)} The Newton-method space-time Jacobian of size $2(N_T+1)N_h^2 \times 2(N_T+1)N_h^2$. 
Each block represents an $N_h^2 \times N_h^2$ spatial operator. 
Note that $D_i$ and $\tilde{E}_i$ are strictly pentadiagonal due to the 2D stencil, while $E_i$ remains diagonal due to local separability.}
\label{fig:sparsity_separableH_2d}
\end{figure}

We next compare the computational costs of the two methods for the $d = 2$ double-well problem.
For the Picard method, the spatial Jacobians are pentadiagonal matrices of size $N_h^2\times N_h^2$ (Figure~\ref{fig:sparsity_separableH_2d}, left), reflecting the 5-point finite-difference stencil of the 2D spatial operators.
For the Newton method, the four block matrices result in a space-time matrix of size $2(N_T+1)N_h^2 \times 2(N_T+1)N_h^2$ (Figure~\ref{fig:sparsity_separableH_2d}, right).
The blocks of $D_i$ in $A_{\cU, \cU}$ and $A_{\cM, \cM}$ are pentadiagonal, as are the blocks of $\tilde{E}_i$ in $A_{\cM, \cU}$.
However, because the Hamiltonian remains locally separable, the blocks in $A_{\cU, \cM}$ remain strictly diagonal.

In two dimensions, each Picard step requires solving a sparse $N_h^2\times N_h^2$ spatial system arising from a five-point stencil. Although the Jacobian has only $O(N_h^2)$ nonzero entries, sparse direct factorization introduces fill-in whose amount depends on the grid structure and the fill-reducing ordering. Table~\ref{tab:scaling_study} reports the resulting fill-in and wall-clock cost for a representative HJB spatial Jacobian factored using SuperLU with COLAMD ordering. The measurements illustrate the increasing cost of direct spatial solves under grid refinement; they are not intended as an empirical determination of an asymptotic complexity rate.

\begin{table}[H]
\centering
\begin{tabular}{l c c c c}
\hline
\textbf{Grid size ($N_h \times N_h$)} & $\mathbf{N_h^2}$ \textbf{(DOFs)} & \textbf{$y$-neighbor offset} & \textbf{Fill-in ratio} & \textbf{Solve time (s)} \\
\hline
$25 \times 25$   & 625    & 25  & 10.9072$\times$ & 0.001256 \\
$50 \times 50$   & 2,500  & 50  & 15.7534$\times$ & 0.023376 \\
$100 \times 100$ & 10,000 & 100 & 23.7258$\times$ & 0.199882 \\
$200 \times 200$ & 40,000 & 200 & 35.9115$\times$ & 2.74088 \\
\hline
\end{tabular}
\caption{Fill-in and direct-solve cost for a representative two-dimensional Picard HJB spatial Jacobian. The matrix was generated from the five-point discretization of the double-well problem with $\nu=0.1$ and $N_T=10$ and then factored using SuperLU with COLAMD ordering. The fill-in ratio is $(\operatorname{nnz}(L)+\operatorname{nnz}(U))/\operatorname{nnz}(J)$, and solve time includes one factorization and one back-substitution.}
\label{tab:scaling_study}
\end{table}

While the shift to 2D significantly increases the cost of the spatial solves,
the relative efficiency of the Picard method becomes clear when contrasted with
the cost of solving the space-time systems arising in the global Newton
method. The space-time coupling in the Newton Jacobian produces sparse linear
systems whose direct factorization rapidly becomes expensive. 

For a representative Newton Jacobian for the two-dimensional double-well
problem with $N_h=50$, $N_T=10$, and $\nu=0.5$, the coupled system has
$2(N_T+1)N_h^2=55{,}000$
DOFs and $455{,}000$ nonzero entries. On an Apple M1 Pro with 16 GB of memory,
SciPy's direct sparse LU factorization \texttt{splu} generated factors
containing $128{,}928{,}646$ nonzeros, corresponding to a $283.36\times$
fill-in ratio, and required $312.43$ seconds.

Before illustrating the hybrid strategy, we characterize each method's reach and wall-clock cost for each $(N_h,\nu)$ run in a sweep over $(N_h,\nu)$ in Tables~\ref{tab:np_2d_nt10} and~\ref{tab:gn_2d_nt10} for the Picard and Newton methods, respectively.
 \begin{table}[h]
  \centering
	  \caption{2D Picard method: wall-clock time (s) for each $(N_h,\nu)$ run, with the damping parameter $\theta$ shown in parentheses. The damping value is retained after each successful viscosity level and decreased in the schedule $\{0.8,0.7,\ldots,0.2,0.1,0.05\}$ when needed. Entries marked with -- indicate failure at every damping value attempted. $N_T = 10$.}
  \label{tab:np_2d_nt10}
  \renewcommand{\arraystretch}{1.15}
  \small
  \begin{tabular}{c|ccccccccc}
  \toprule
  $N_h \backslash \nu$ & $1.0$ & $0.3$ & $0.1$ & $0.03$ & $0.01$ & $0.003$ & $0.001$ & $3\!\times\!10^{-4}$ & $10^{-4}$ \\
  \midrule
   $16$  & $0.5\,(0.8)$ & $0.5\,(0.8)$ & $0.5\,(0.8)$ & $4.9\,(0.5)$ & $4.7\,(0.3)$ & $2.3\,(0.3)$ & $2.3\,(0.3)$ & $2.4\,(0.3)$ & $2.3\,(0.3)$ \\
    $32$  & $1.2\,(0.8)$ & $1.3\,(0.8)$ & $1.3\,(0.8)$ & $8.6\,(0.5)$ & $12.6\,(0.3)$ & $14.3\,(0.2)$ & -- & -- & -- \\
    $64$  & $5.2\,(0.8)$ & $5.6\,(0.8)$ & $6.4\,(0.8)$ & $39.9\,(0.5)$ & $51.1\,(0.3)$ & -- & -- & -- & -- \\
    $128$ & $39.9\,(0.8)$ & $43.7\,(0.8)$ & $76.7\,(0.8)$ & $616.1\,(0.5)$ & $413.9\,(0.3)$ & -- & -- & -- & -- \\
  \bottomrule
  \end{tabular}
\end{table}

  \begin{table}[h]
  \centering
	  \caption{2D Newton method: cumulative wall-clock time to reach each $\nu$ via viscosity continuation from $\nu = 1.0$. $N_T = 10$. }
  \label{tab:gn_2d_nt10}
  \renewcommand{\arraystretch}{1.15}
  \begin{tabular}{c|ccccccccc}
  \toprule
  $N_h \backslash \nu$ & $1.0$ & $0.3$ & $0.1$ & $0.03$ & $0.01$ & $0.003$ & $0.001$ & $3\!\times\!10^{-4}$ & $10^{-4}$ \\
  \midrule
  $16$  & $1.4$ & $2.3$ & $3.2$ & $3.8$ & $4.4$ & $5.0$ & $5.4$ & $5.9$ & $6.3$ \\
  $32$  & $34.9$ & $57.5$ & $82.8$ & $106.2$ & $127.8$ & $148.3$ & $163.8$ & $174.2$ & $184.5$ \\
  $64$  & $632.8$ & $988.4$ & $1{,}344.5$ & $1{,}713.1$ & $2{,}079.9$ & $2{,}449.3$ & $2{,}730.9$ & $3{,}009.1$ & $3{,}289.2$ \\
  $128$ & $7{,}032.8$ & $12{,}061.5$ & $18{,}171.6$ & $23{,}203.2$ & $28{,}672.5$ & $33{,}524.0$ & $38{,}452.6$ & $42{,}384.0$ & $46{,}410.5$ \\
  \bottomrule
  \end{tabular}
  \end{table}

Table~\ref{tab:np_2d_nt10} shows that Picard's reachable viscosity depends on spatial resolution.
On the coarsest grid, $N_h=16$, the method reaches $\nu=10^{-4}$, consistent with the additional stabilization supplied by numerical diffusion in the upwind discretization.
This advantage diminishes under refinement: the smallest successful viscosity is $3\times10^{-3}$ for $N_h=32$, and $10^{-2}$ for both $N_h=64$ and $128$.
Thus, on the finer grids tested here, Picard exhibits an apparent viscosity floor under the prescribed damping schedule and stopping criteria.
This motivates the hybrid strategy: use Picard above that floor, then switch to Newton continuation to reach smaller viscosities.

For the Newton method, the first cell in each row ($\nu = 1$) starts from the naive guess $U \equiv u_T$, $M \equiv m_0$ and requires $7$ Newton updates (global linear solves) to satisfy the convergence tolerance. Each subsequent solve starts from the preceding converged solution and requires $2$--$5$ updates.

We now illustrate the hybrid strategy of Section~\ref{sec:numericalmethods:hybrid} with the target $(N_h, \nu_{\text{target}}) = (128, 10^{-4})$. 
Pure Picard fails at this target; pure Newton reaches the target via continuation at a cumulative cost of $46{,}410$\,s $\approx 12.9$ hours (Table~\ref{tab:gn_2d_nt10}). 
The hybrid traverses a single row of the $(N_h, \nu)$ grid, summarized in Table~\ref{tab:hybrid_phases}:

  \begin{table}[H]
    \centering
    \caption{Hybrid run at the target row $N_h = 128$, $\nu_{\text{target}}=10^{-4}$,
    resolved by phase. \colorbox{Green}{Phase 1}: Picard descent, one row per
    viscosity reached. \colorbox{Salmon}{Phase 2}: the Picard failure at the next
    viscosity, which triggers the handoff and identifies
    $\nu_{\text{limit}}$. \colorbox{SkyBlue}{Phase 3}: Newton, comprising
    the re-solve at $\nu_{\text{limit}}$ followed by continuation to
    $\nu_{\text{target}}$. Iteration counts are Picard iterations in Phases 1--2
    and Newton updates in Phase 3. The damping parameter $\theta$ used for each Picard run is given in parentheses beside the method.}
    \label{tab:hybrid_phases}
    \begin{tabular}{c c l r r}
      \toprule
      \textbf{Phase} & $\boldsymbol{\nu}$ & \textbf{Method} & \textbf{Iterations}
        & \textbf{Cumulative (s)} \\
      \midrule
      \cellcolor{Green}       & $1.0$              & Picard $(\theta = 0.8)$ & $11$ & $39.1$ \\
      \cellcolor{Green}       & $0.3$              & Picard $(\theta = 0.8)$ & $11$ & $82.4$ \\
      \cellcolor{Green}1      & $0.1$              & Picard $(\theta = 0.8)$ & $19$ & $159.8$ \\
      \cellcolor{Green}       & $0.03$             & Picard $(\theta = 0.5)$ & $102$ & $785.2$ \\
      \cellcolor{Green}       & $0.01$             & Picard $(\theta = 0.3)$ & $51$ & $1{,}221.5$ \\
      \midrule
      \cellcolor{Salmon}2     & $0.003$ & Picard $(\theta = 0.05)$ & $11$\rlap{$^{\dagger}$} & $1{,}644.3$ \\
      \midrule
      \cellcolor{SkyBlue}     & $0.01$              & Newton & $1$ & $2{,}809.2$ \\
      \cellcolor{SkyBlue}     & $0.003$             & Newton & $4$ & $8{,}084.8$ \\
      \cellcolor{SkyBlue}3    & $0.001$             & Newton & $4$ & $12{,}922.8$ \\
      \cellcolor{SkyBlue}     & $3\!\times\!10^{-4}$ & Newton & $3$ & $16{,}237.5$ \\
      \cellcolor{SkyBlue}     & $10^{-4}$           & Newton & $3$ & $19{,}563.5$ \\
      \bottomrule
  \end{tabular}

    \vspace{2pt}
    {\footnotesize $^{\dagger}$Did not converge: at $\theta = 0.05$, no
    improvement of the best residual occurred for ten consecutive iterations;
    the stall criterion therefore terminated the run at iteration $11$, before
    the $200$-iteration cap.}
  \end{table}

The speedup arises from the work avoided by the hybrid strategy. 
Pure Newton spends roughly $62\%$ of its total cost on the high-$\nu$ continuation steps reaching $\nu_{\text{limit}}$ --- steps the hybrid strategy replaces with the much cheaper Picard sweep. 
The hybrid reaches the target where pure Picard fails by switching to Newton at the Picard method's lower-viscosity limit, while reducing the cumulative runtime relative to pure Newton by replacing its expensive high-viscosity continuation steps with cheaper Picard iterations.

\section{Conclusion and future work}
\label{sec:conc}

We have developed and tested a hybrid strategy that uses Picard iterations to accelerate Newton continuation for finite-difference MFG systems. The computational savings come from replacing expensive high-viscosity Newton solves with inexpensive Picard solves during the initial viscosity descent. The method then transfers the last converged Picard solution to Newton to reach a smaller target viscosity on the same discretization. The component-method studies provide the basis for this combination by identifying where Picard offers a computational advantage and where Newton's greater robustness becomes necessary.

The numerical examples also identify circumstances in which a hybrid strategy offers little benefit. For nonlocal congestion, dense spatial coupling makes assembly of the full Newton Jacobian impractical at the grid sizes considered, while Picard remains viable through FFT-based evaluation of the nonlocal terms. In the traffic-light problem, strong damping reduces Picard's computational advantage: the recorded Newton run reaches a lower viscosity in less time than the recorded Picard run. These examples clarify when the Picard-accelerated hybrid method is appropriate: when Picard is inexpensive over part of the parameter range and Newton is needed to reach the target regime.

These computational gains are established for the common direct sparse backend used in the reported experiments. An important extension is to assess the hybrid method with alternative linear solvers, including the block strategies of Achdou and Perez~\cite{achdou2012preconditioning} for Newton and incomplete-LU preconditioning for Picard's spatial systems. Such a study should distinguish changes in the number of linear solves from changes in their individual cost, since improving the linear algebra can change the relative advantage of each phase and the appropriate handoff point.

Further work should investigate a handoff criterion based on computational cost, rather than waiting for the Picard damping schedule to be exhausted. An intermediate Picard/Newton variant also merits study: retaining the outer HJB--FP decoupling while solving the HJB equation globally in time would place the method between timestep-by-timestep Picard sweeps and the fully coupled Newton solve. Its convergence properties and computational cost remain to be assessed.

Finally, viscosity continuation should be studied alongside line-search globalization for Newton's method. Continuation supplies initial guesses by changing the viscosity, whereas an Armijo line search, as used by Carlini and Zorkot~\cite{Carlini2027Newton}, controls Newton steps at fixed parameters. Their combination could broaden the regimes in which Picard-based acceleration is effective. Establishing convergence criteria and selecting the handoff and continuation schedules from computable diagnostics are natural next steps toward a more systematic hybrid method.

The code used to generate the figures and timings is available at \href{https://github.com/andrewjshi/2026-pap-Newton}{https://github.com/andrewjshi/2026-pap-Newton}.

\bibliographystyle{elsarticle-num}
\bibliography{achdou, pap}

\appendix
\section{Pseudocode}

\begin{algorithm}
\caption{Newton method}
\label{alg:global-newton}
\begin{algorithmic}[1]
\State \textbf{Initialize:} Guess $\mathcal{W}^{(0)} = (U^{(0)}, M^{(0)})^{\top}$, set tolerance $\epsilon_N$, $k \gets 0$
\While{not converged}
    \State \textbf{1. Evaluate residuals:}
    \State Compute $\varphi(\mathcal{W}^{(k)}) = \begin{pmatrix} \varphi_{U}(U^{(k)}, M^{(k)}) \\ \varphi_{M}(U^{(k)}, M^{(k)}) \end{pmatrix}$ using \eqref{eq:MFG-finitedifference-fb:HJB} and \eqref{eq:MFG-finitedifference-fb:FP} 
    
    \State \textbf{2. Construct the Jacobian:}
    \State Assemble $J_{\varphi}(\mathcal{W}^{(k)}) = \begin{pmatrix} A_{\cU,\cU} & A_{\cU,\cM} \\ A_{\cM,\cU} & A_{\cM,\cM} \end{pmatrix}$ 
    
    \State \textbf{3. Solve the linear system for the update:}
    \State Solve $J_{\varphi}(\mathcal{W}^{(k)}) \delta \mathcal{W}^{(k)} = -\varphi(\mathcal{W}^{(k)})$ for $\delta \mathcal{W}^{(k)} = (\tilde{U}, \tilde{M})^{\top}$ 
    
    \State \textbf{4. Update the solution:}
    \State $\mathcal{W}^{(k+1)} \gets \mathcal{W}^{(k)} + \delta \mathcal{W}^{(k)}$ 
    
    \State \textbf{5. Check convergence:}
    \If{$\|\varphi(\mathcal{W}^{(k+1)})\| < \epsilon_N$}
        \State \textbf{return} $(U^{(k+1)}, M^{(k+1)})$
    \EndIf
    \State $k \gets k + 1$
\EndWhile
\end{algorithmic}
\end{algorithm}

\begin{algorithm}
\caption{Picard method: damped Picard with backward HJB Newton and forward FP}
\label{alg:picard}
\begin{algorithmic}[1]
\State \textbf{Initialize:} Set $U^{(0)} \equiv 0$ and $M^{(0)}$ from $m_0$, set Picard tolerance $\epsilon_{P}$, Newton tolerance $\epsilon_{N}$, damping factor $\theta \in (0, 1]$, $k \gets 0$
\While{not converged}
\State \textbf{1. Solve nonlinear HJB for $\hat{U}$ (backward-marching Newton):} 
    \State Set terminal condition: $\hat{U}^{N_T} = g(x, M^{(k), N_T})$ 
    \For{$n = N_T - 1$ down to $0$} 
        \State Find $\hat{U}^n$ by solving the localized nonlinear HJB equation \eqref{eq:MFG-finitedifference-fb:HJB}.
        \State \textbf{Use next time step as initial guess:} Let $V^{(0)} \gets \hat{U}^{n+1}$ 
        \For{$\ell_{\mathrm{N}} = 0, 1, \dots$} 
            \State Compute local residual $R_{\ell_{\mathrm{N}}} = \text{Res}_n(V^{(\ell_{\mathrm{N}})}, \hat{U}^{n+1}, M^{(k), n+1})$
            \If{$\|R_{\ell_{\mathrm{N}}}\| < \epsilon_N$} \State \textbf{break} \EndIf
            \State Solve local Jacobian $J_n(V^{(\ell_{\mathrm{N}})}) \delta u = -R_{\ell_{\mathrm{N}}}$
            \State $V^{(\ell_{\mathrm{N}}+1)} \gets V^{(\ell_{\mathrm{N}})} + \delta u$
        \EndFor
        \State $\hat{U}^n \gets V^{(\mathrm{final})}$
    \EndFor
    \State \textbf{Damped update for $U$:} $U^{(k+1)} \gets (1-\theta)U^{(k)} + \theta \hat{U}$

    \State \textbf{2. Solve linear Fokker--Planck for $\hat{M}$ (forward step):} 
    \State Given $U^{(k+1)}$, solve the FP equation for $\hat M$, evaluating density-dependent coefficients in $\partial\tilde H/\partial q_\ell$ at $M^{(k)}$.
    \State March forward in time from $M^0 = m_0$:
    \State Let $\mathcal{B}^{(k)}(U,\cdot)$ denote the transport operator
    \eqref{eq:transport-operator} with the density-dependent coefficients
    $\partial\tilde H/\partial q_\ell$ evaluated at $M^{(k)}$.
    \State $\frac{\hat{M}_{i,j}^{n+1}-\hat{M}_{i,j}^n}{\Delta t}
    -\nu(\Delta_h\hat{M}^{n+1})_{i,j}
    -\mathcal{B}^{(k)}_{i,j}(U^{(k+1),n},\hat{M}^{n+1})=0$
    \State \textbf{Damped update for $M$:} $M^{(k+1)} \gets (1-\theta)M^{(k)} + \theta \hat{M}$
    
    \State \textbf{3. Outer convergence check:} 
    \If{$k\geq 1$ \textbf{and}
    $\frac{\|U^{(k+1)}-U^{(k)}\|}{\|U^{(k)}\|}<\epsilon_P$
    \textbf{and}
    $\frac{\|M^{(k+1)}-M^{(k)}\|}{\|M^{(k)}\|}<\epsilon_P$}
        \State \textbf{return} $(U^{(k+1)}, M^{(k+1)})$
    \EndIf
    \State $k \gets k + 1$
\EndWhile
\end{algorithmic}
\end{algorithm}

\begin{algorithm}
\caption{Newton method with viscosity continuation}
\label{alg:newton_continuation}
\begin{algorithmic}[1]
\State \textbf{Initialize:} 
\State Choose starting viscosity $\nu_{\text{curr}}$ (typically $0.5$ or $1.0$ per \cite{achdou2012preconditioning})
\State Choose target viscosity $\nu_{\text{target}}$
\State Set reduction factor $\tau \in (0, 1)$ (e.g., $0.5$)
\State Set backtracking factor $\gamma \in (0, 1)$ (e.g., $0.5$)
\State Initial guess: $\mathcal{W} \leftarrow (U_{T}, M_{0})$ for all times

\State \textbf{First solve:}
\State $\mathcal{W} \leftarrow \text{NewtonSolver}(\nu_{\text{curr}}, \mathcal{W})$ \Comment{Run Algorithm~\ref{alg:global-newton}}
\If{diverged} \State \textbf{abort} \EndIf

\While{$\nu_{\text{curr}} > \nu_{\text{target}}$}
    \State $\nu_{\text{next}} \leftarrow \max(\nu_{\text{target}}, \tau \nu_{\text{curr}})$
    \State success $\leftarrow$ \textbf{false}
    
    \While{\textbf{not} success}
        \State $\mathcal{W}_{trial} \leftarrow \text{NewtonSolver}(\nu_{\text{next}}, \mathcal{W})$ \Comment{Use previous solution as guess}
        
        \If{converged}
            \State $\mathcal{W} \leftarrow \mathcal{W}_{trial}$
            \State $\nu_{\text{curr}} \leftarrow \nu_{\text{next}}$
            \State success $\leftarrow$ \textbf{true}
        \Else
            \State \textbf{Backtrack:} Reduce step size
            \State $\nu_{\text{next}} \leftarrow \nu_{\text{curr}} - \gamma (\nu_{\text{curr}} - \nu_{\text{next}})$
            \If{$|\nu_{\text{curr}} - \nu_{\text{next}}| < \text{tol}$} 
                \State \textbf{abort} (step size too small)
            \EndIf
        \EndIf
    \EndWhile
\EndWhile

\State \textbf{return} $\mathcal{W}$
\end{algorithmic}
\end{algorithm}

\begin{algorithm}
    \caption{Hybrid strategy: combined Picard and Newton methods}
    \label{alg:hybrid_strategy}
    \begin{algorithmic}[1]
     \State \textbf{Inputs:} Target discretization $(N_h,N_T)$, viscosity schedule
     $\nu_{\text{start}} = \nu_0>\nu_1>\cdots>\nu_L=\nu_{\text{target}}$
    \State \textbf{Initialize:}
    \State $\nu \leftarrow \nu_{\text{start}}$; $\nu_{\text{limit}} \leftarrow
  \nu_{\text{start}}$
    \State $\theta \leftarrow 0.8$; $\theta_{\text{floor}} \leftarrow 0.05$;
  $\theta_{\text{step}} \leftarrow 0.1$
    \State $\mathcal{W} \leftarrow \varnothing$ \Comment{No converged Picard state yet}
    \Statex \textbf{Phase 1: Picard descent in $\nu$ at the target discretization}
    \While{$\nu \geq \nu_{\text{target}}$}
        \State $\mathcal{W}_{\mathrm{init}} \leftarrow \operatorname{InitPicard}(\nu)$
        \Comment{Zero value; density from pure diffusion}
        \State $\mathcal{W}_{\mathrm{trial}} \leftarrow \text{PicardSolver}(\nu, (N_h, N_T),
  \mathcal{W}_{\mathrm{init}}, \theta)$
        \If{converged}
            \State $\mathcal{W} \leftarrow \mathcal{W}_{\mathrm{trial}}$
            \State $\nu_{\text{limit}} \leftarrow \nu$
            \Comment{Update the last successful viscosity level}
            \If{$\nu=\nu_{\text{target}}$}
                \State \textbf{return} $\mathcal{W}$
                \Comment{Target reached with the Picard method}
            \EndIf
            \State $\nu \leftarrow \text{NextViscosity}(\nu)$
            \Comment{Move to the next viscosity in the descent}
  \Else \Comment{Discard the failed trial; retain the last converged state}
      \If{$\theta \le \theta_{\text{floor}}$}
          \Statex \textbf{Phase 2: Picard failure and handoff}
          \State \textbf{break} \Comment{Retain the last converged state at $\nu_{\text{limit}}$}
      \EndIf
      \State $\theta \leftarrow \max(\theta_{\text{floor}}, \theta - \theta_{\text{step}})$
  \Comment{Sticky decrement; retry same $\nu$}
  \EndIf
    \EndWhile
    \Statex \textbf{Phase 3: Newton re-solve and continuation to $\nu_{\text{target}}$}
    \State $\mathcal{W} \leftarrow \text{NewtonSolver}(\nu_{\text{limit}}, (N_h, N_T),
  \mathcal{W})$ \Comment{Reduce the coupled residual}
    \State $\nu \leftarrow \nu_{\text{limit}}$ \Comment{Resume the descent from the handoff
  viscosity}
    \While{$\nu > \nu_{\text{target}}$}
        \State $\nu \leftarrow \text{NextViscosity}(\nu)$ \Comment{Move to next $\nu$ in the
  descent}
        \State $\mathcal{W} \leftarrow \text{NewtonSolver}(\nu, (N_h, N_T), \mathcal{W})$
    \EndWhile
    \State \textbf{return} $\mathcal{W}$
    \end{algorithmic}
    \end{algorithm}    
    
\end{document}